\documentclass[amsfonts,12pt]{amsart}

\usepackage{cite}
\usepackage{epsfig}
\usepackage{mathtools}
\usepackage[subnum]{cases}

\usepackage{amsmath,amssymb}
\usepackage[bookmarks=false,
bookmarksnumbered=true, breaklinks=true,
pdfstartview=FitH, hyperfigures=false,
plainpages=false, naturalnames=true,
colorlinks=true,pagebackref=true,
pdfpagelabels]{hyperref}

\allowdisplaybreaks
\newtheorem{theorem}{Theorem}[section]

\newtheorem{lemma}[theorem]{Lemma}
\newtheorem{corollary}[theorem]{Corollary}

\newtheorem{proposition}[theorem]{Proposition}

\begin{document}
	\title{The $L_p$ Minkowski problems for affine dual quermassintegrals  }
	\author[Youjiang Lin, Yuchi Wu]
	{Youjiang Lin, Yuchi Wu}
	\address{School of Mathematical Sciences, Hebei Normal University, Shijiazhuang, Hebei,
050024, PR China} \email{yjlin@hebtu.edu.cn}
	\address{School of Mathematical Sciences, Key Laboratory of MEA(Ministry of Education) \& Shanghai Key Laboratory of PMMP, East China Normal University, Shanghai, 200241, PR China} \email{ycwu@math.ecnu.edu.cn}

	\thanks{ There is no conflict of interest.
The first author is supported by National Natural Science Foundation of China NSFC 12371137 and NSFC 11971080. The second author is supported by Science and Technology Commission of Shanghai Municipality 22DZ2229014 and Youth Fund of the National Natural Science Foundation of China NSFC 12401067}
	\subjclass[2020]{Primary: 52A20, 52A30, 52A40} \keywords{Affine convex geometry; $L_p$ Minkowski problem; Dual affine quermassintegral; Affine
dual curvature measure; Radon transform\; Grassmannian}
	\maketitle
	\pagestyle{myheadings}
	\markboth{YOUJIANG LIN AND YUCHI WU}{The $L_p$ Minkowski problems on ADQ}
	
	\begin{abstract}
 In this paper,  {we extend
the affine dual curvature measures to the $L_p$ setting and solve the existence part of the corresponding Minkowski problem  for non-symmetric discrete measures when $p>1$ and for symmetric measures when $p\geq0$.} When $p=0$, the $L_0$ Minkowski problem is the affine dual Minkowski problem, which is introduced and solved by Cai-Leng-Wu-Xi in \cite{CLWX}.
	\end{abstract}

%
%
%
%
%
%
%
%
%
%
	
	\section{Introduction}

 {A Minkowski problem concerns the identification of criteria that determine whether a specified measure can be expressed as the geometric measure associated with a convex body. In essence, it investigates the conditions under which such a measure originates from a convex body}. The
Minkowski problem is one of the cornerstones of the Brunn-Minkowski theory of convex bodies. And, if the convex body exists, to what extent is it unique? This is  the uniqueness question for the Minkowski problem. There are also regularity questions as well as stability questions for the Minkowski problem. Many very important research achievements on Minkowski problems have been made so far, see e.g., \cite{GMZ06,HXZ21,Lutwak93,LXYZ23,Xiong10,LYZ24,Zhu14,Zhu15}.

$L_p$-Minkowski problem is to find necessary and sufficient conditions on a finite Borel measure $\mu$ on the unit sphere $S^{n-1}$ so that $\mu$ is the $L_p$-surface area measure of a convex body in $\mathbb{R}^n$. The case $p=1$ of the $L_p$-Minkowski problem is of course the classical Minkowski problem, which was solved by Minkowski, Alexandrov, and Fenchel and Jessen (see Schneider \cite{Schneider14} for references). In addition to the surface area measure of a convex body, another fundamental measure associated with a convex body $K$ in $\mathbb{R}^n$ that contains the origin in its interior is the cone-volume measure. A very important property of the cone-volume measure is its $SL(n)$ invariance, or
simply called affine invariance. The $L_p$-surface area measures for $p\neq 0$ are all $\operatorname{SO}(n)$ invariant. The  Minkowski
problem for the cone-volume measure is called the logarithmic Minkowski problem. It asks for necessary and sufficient conditions for a given measure on the unit sphere to be the cone-volume measure of a convex body. B\"or\"oczky-Lutwak-Yang-Zhang \cite{BLYZ13} solved the existence part of the logarithmic Minkowski problem  for the case of even measures within the
class of origin-symmetric convex bodies.  On more the known solutions of the $L_p$ Minkowski problems, see e.g., \cite{CLZ17,CW06,GXZ,HS09,HLYZ18,HYZ25,Lutwak93,WXL19,Zhu15}.

Lutwak \cite{Lutwak75,Lutwak88} built the dual Brunn-Minkowski theory as a `dual' counterpart of the classical Brunn-Minkowski theory. In the dual Brunn-Minkowski theory,  Huang-Lutwak-Yang-Zhang \cite{HLYZ16} discovered a new family of geometric measures called dual curvature measures. The dual Minkowski problem is a class of  Minkowski problems  associated with
the dual curvature measures, and is a major problem in the dual Brunn-Minkowski theory. Existence of solutions for the
dual Minkowski problem for even data within the class of origin-symmetric convex bodies was established in \cite{HLYZ16}. Recently existence for the critical cases of the even dual Minkowski problem were established in \cite{BHP18} and \cite{Zhao17}, and a complete solution to the dual
Minkowski problem with negative indices was given in \cite{Zhao18}. On more the known solutions of the $L_p$ dual Minkowski problems, see e.g., \cite{BLYZZ19,CLZ17,CW06,GHXYW20,HP18,HZ18,LSW18,Zhang94,Zhao17}.

The affine geometry, and centro-affine geometry, of convex bodies studies affine and centro-affine invariants, i.e., geometric functionals that are invariant under affine and and centro-affine transformations. There are notions in affine convex geometry such as affine surface area, the Mahler volume, the projection, intersection, and centroid operators
that are fundamental in this fruitful theory. In the 1980’s the classical affine surface area
of smooth convex bodies from affine differential geometry was finally extended to all convex
bodies, see e.g. \cite{Lutwak91,Meyer00,SW04,WY10} and the book of Leichtweiss \cite{Lei98}. Various new affine invariants were studied during the past half century, see e.g., \cite{GG99,Gardner9401,Gardner9402,Haberl,HP14,Kol00,Lin17,LX21,LW22,Ludwig03,Ludwig06,Ludwig10,LYZ0001,Lutwak96,LYZ00,LYZ02,Zhang96,Zhang99}.


 Affine quermassintegrals \cite{Lutwak90} and affine  dual  quermassintegrals were both proposed by Lutwak for a convex body $K\in \mathcal{K}^n_{(o)}$, see \cite{Schneider14}. They were defined by letting $\Phi_0(K)=\widetilde{\Phi}_0(K)=V(K)$, and $\Phi_n(K)=\widetilde{\Phi}_n(K)=\omega_n$, while for $m=1, \ldots, n-1$,
$$
\Phi_{n-m}(K):=\frac{\omega_n}{\omega_m}\left(\int_{G(n, m)} \operatorname{vol}_m(K \mid \xi)^{-n} d \nu_m(\xi)\right)^{-1 / n}
$$
and
\begin{equation}\label{Ddacm}
\widetilde{\Phi}_{n-m}(K):=\frac{\omega_n}{\omega_m}\left(\int_{G(n, m)} \operatorname{vol}_m(K \cap \xi)^n d \nu_m(\xi)\right)^{1 / n},
\end{equation}
where $G(n, m)$ denotes the Grassmannian of $m$ dimensional subspaces of $\mathbb{R}^n$, while $\nu_m$ denotes the Haar measure (the unique rotationally invariant probability measure) on $G(n, m)$, while $\operatorname{vol}_m$ denotes $m$ dimensional volume, and $\omega_m$ denotes the volume of the unit ball in $\mathbb{R}^m$.
The normalization is to preserve the geometric meaning of these quantities.
It was shown by Grinberg \cite{Grinberg91} that $\Phi_{n-m}$ and $\widetilde{\Phi}_{n-m}$ are, as their names suggests, $SL(n)$ invariant. Moreover, $\Phi_{n-m}$ is translation invariant as well. Isoperimetric inequalities for these affine invariants are stronger than their classical counterparts. For the affine quermassintegral $\Phi_{n-m}(K)$, it had been a long-standing conjecture of Lutwak that the lower bound is achieved by ellipsoids.  Lutwak's conjecture was verified  recently by Milman-Yehudayoff \cite{MY24} after more than three decades of efforts.

Let $m \in\{1, \ldots, n-1\}$. The {\it $m$-dimensional Radon transform} $\mathcal{R}_m$ is a linear operator from $C\left(S^{n-1}\right)$ into $C(G(n, m))$ given by
\begin{equation*}\label{Rm}
\mathcal{R}_m f(\xi)=\int_{S^{n-1} \cap \xi} f(w) d w.
\end{equation*}
Here $d w$ is used as an abbreviation for $d \mathcal{H}^{m-1}(w)$ on $S^{n-1} \cap \xi$. The {\it $m$ dimensional dual Radon transform} $\mathcal{R}_m^*: C(G(n, m)) \rightarrow C\left(S^{n-1}\right)$ is defined by
\begin{equation*}\label{Rmstar}
\mathcal{R}_m^* F(u)=\frac{m \omega_m}{n \omega_n} \int_{G_u(n-1, m-1)} F(\operatorname{span}\{u, \zeta\}) d \nu_{m-1}(\zeta)
\end{equation*}
where $G_u(n-1, m-1)$ denotes the Grassmannian of $m-1$ dimensional subspaces of $u^{\perp}$ and $\nu_{m-1}$ is the Haar measure on $G_u(n-1, m-1)$, and as always, we use ``Haar measure" to mean ``Haar probability measure".

In \cite{CLWX}, Cai-Leng-Wu-Xi constructed  the affine dual curvature measure: Suppose $m \in\{1, \ldots, n-1\}$ and $K\in \mathcal{K}^n_{(o)}$, the affine dual curvature measure of $K$ is defined by
\begin{equation}\label{Dadcm}
\widetilde{C}_m^{\mathrm{a}}(K, \eta)=\int_{\boldsymbol{\alpha}_K^*(\eta)} \rho_K(u)^m \mathcal{R}_m^*\left(\operatorname{vol}_m(K \cap \cdot)^{n-1}\right)(u)d u
\end{equation}
for each Borel set $\eta \subset S^{n-1}$, where $\boldsymbol{\alpha}_K^*(\eta)$ denotes the reverse radial Gauss image of $\eta$ (see Section \ref{s2.2} for details), $\mathcal{R}_m^*$ denotes $m$ dimensional dual
Radon transform (\ref{Rmstar}). The affine dual curvature measures of convex bodies are affine contravariant (see \cite{BLYZ15}).  More precisely, for each $\varphi \in SL(n)$,
$$
\widetilde{C}_m^{\mathrm{a}}(\varphi K, \cdot)=\varphi^{-\mathrm{t}} \widetilde{C}_m^{\mathrm{a}}(K, \cdot).
$$
Moreover, Cai-Leng-Wu-Xi \cite{CLWX} proved the variation formula for affine dual curvature measures, which implies that affine dual curvature measure is obtained from taking the derivative of affine dual quermassintegral.

The {\it affine dual quermassintegrals} defined in (\ref{Ddacm}) can been extended to the compact convex sets that may not contain origin in their interiors, i.e., for $K\in\mathcal{K}_o^n$,
\begin{equation*}
\widetilde{\Phi}_{n-m}(K):=\frac{\omega_n}{\omega_m}\left(\int_{G(n, m)} \operatorname{vol}_m(K \cap \xi)^n d \nu_m(\xi)\right)^{1 / n}.
\end{equation*}
And let
\begin{equation*}\label{2.2}
\widetilde{\Psi}_{m}(K)=\int_{G(n, m)}\operatorname{vol}_m(K \cap \xi)^n d \nu_m(\xi).
\end{equation*}
It is obvious that $\widetilde{\Phi}_{n-m}(K)=0$ and $\widetilde{\Psi}_{m}(K)=0$ when $\operatorname{dim}K\leq n-1$. Next, for abbreviation, we write ``ADQ" denote ``affine dual quermassintegrals".

Moreover, we define $\widetilde{C}_{p,m}^a(K, \cdot)$ to be the {\it $L_p$ curvature measure of ADQ}: Let $p\in\mathbb{R}$, $n\geq 2$, $m \in\{1, \ldots, n-1\}$, and $K\in \mathcal{K}^n_{o}$ with $\operatorname{int} K \neq \emptyset$. Then
\begin{equation*}
\widetilde{C}_{p,m}^a(K, \eta):=\int_{\boldsymbol{\alpha}_K^*(\eta)\cap\left(\operatorname{int} N(K, o)^*\right)} \left(u\cdot\nu_K\left(\rho_K(u)u\right)\right)^{-p}\rho_K(u)^{m-p} \mathcal{R}_m^*\left(\operatorname{vol}_m(K \cap \cdot)^{n-1}\right)(u) d u
\end{equation*}
for each Borel set $\eta \subset S^{n-1}$, where $N(K, o)^*$ denotes the dual of  the normal cone $N(K,o)$ (see Section \ref{s2.2} for details).

By the above definition of the  $L_p$ curvature measure of ADQ, for $K\in \mathcal{K}^n_{(o)}$ and each continuous $g:S^{n-1}\rightarrow \mathbb{R}$,
\begin{equation*}\label{E8}
\int_{S^{n-1}}g(u)d\widetilde{C}_{p,m}^a(K, u)=\int_{S^{n-1}}g\left(\alpha_K(u)\right) h_K\left(\alpha_K(u)\right)^{-p} \rho_K(u)^m \mathcal{R}_m^*\left(\operatorname{vol}_m(K \cap \cdot)^{n-1}\right)(u) d u,
\end{equation*}
where $\alpha_K$ is the radial Gauss map (see Section \ref{s2.2} for details).

For $K\in \mathcal{K}^n_{(o)}$, the $L_p$ curvature measure $\widetilde{C}_m^{\mathrm{a}}(K, \cdot)$ has a geometric interpretation by way of Zhang intersection bodies. In fact, for $m=n-1$, recalling (\ref{Rmstar}), it can be shown that, for $u \in S^{n-1}$,
$$
\mathcal{R}_m^*\left(\operatorname{vol}_m(K \cap \cdot)\right)(u)=\frac{n-1}{n \omega_n} \rho_{\mathrm{I}^2 K}(u),
$$
where $\mathrm{I}^2 K=\mathrm{I}(\mathrm{I} K)$ is the intersection body of $\mathrm{I} K$ (see Section \ref{s2.3} for details). Therefore
\begin{equation*}
\widetilde{C}_{p,n-1}^{\mathrm{a}}(K, \eta)=\frac{n-1}{n \omega_n} \int_{\boldsymbol{\alpha}_K^*(\eta)} h_K\left(\alpha_K(u)\right)^{-p}\rho_K(u)^{n-1} \rho_{\mathrm{I}^2 K}(u) d u.
\end{equation*}

 For generic $m$ and $\xi \in G(n, m)$, by the definition of bi-dual $m$-intersection body $\mathbb{I}_mK$ of $K\in\mathcal{K}_{(o)}^n$ (see (\ref{Dbib}) in Section \ref{s2.3}), we have
\begin{equation}\label{mathI}
\mathcal{R}_m^*\left(\operatorname{vol}_m(K \cap \cdot)^{n-1}\right)(u)= \rho_{\mathbb{I}_m K}(u)^{n-m}.
\end{equation}

By (\ref{mathI}), we can reformulate (\ref{E8}) as
\begin{equation*}
\int_{S^{n-1}}g(u)\widetilde{C}_{p,m}^{\mathrm{a}}(K, u)=\int_{S^{n-1}}g\left(\alpha_K(u)\right) h_K\left(\alpha_K(u)\right)^{-p} \rho_K(u)^m \rho_{\mathbb{I}_m K}(u)^{n-m} d u.
\end{equation*}
By the notion of mixed dual curvature measures introduced by Lutwak-Yang-Zhang \cite{LYZ18} (see Section \ref{s2.4} for details), we can write our new measures $\widetilde{C}_{p,m}^{\mathrm{a}}(K, \cdot)$ as mixed dual curvature measures:
\begin{equation*}
\widetilde{C}_{p,m}^{\mathrm{a}}(K, \cdot)=\widetilde{C}_{p,m}\left(K, \mathbb{I}_m K, \cdot\right).
\end{equation*}

 Note that just as the $L_p$ surface area measure does not have affine invariance for $p\neq 0$, the $L_p$ curvature measure of ADQ also does not have affine invariance for $p\neq0$. When $p=0$, the $L_0$ curvature measure of ADQ  is the  affine dual curvature measure $\widetilde{C}_{m}^a(K, \cdot)$ defined in  (\ref{Dadcm}). \\

\noindent{\bf $L_p$ Minkowski problem for ADQ.}
Suppose $m=1,\dots,n-1$ is fixed and $p\in\mathbb{R}$. Find necessary and sufficient conditions that a finite Borel $\mu$ on $S^{n-1}$ must satisfy for there to exist a convex body $K\in \mathcal{K}^n_{o}$ with $\operatorname{int}K\neq \emptyset$ such that
$$
\mu=\widetilde{C}_{p,m}^{\mathrm{a}}(K, \cdot) .
$$
 {When the measure} $\mu$ has a density function $g: S^{n-1} \rightarrow \mathbb{R}$, the corresponding partial differential equation is a Monge-Ampère type equation on $S^{n-1}$:
\begin{eqnarray}
\mathcal{R}_m^* \circ \mathcal{R}_m^{(n-1)} \circ \mathcal{D}^{(-m)}(h)\left(\frac{\nabla h(v)+h(v) v}{|\nabla h(v)+h(v) v|}\right)\cdot|\nabla h(v)+h(v) v|^{m-n} h(v)^{1-p} \operatorname{det}\left(\nabla^2 h(v)+h(v) I\right)\nonumber\\=g(v),\nonumber
\end{eqnarray}
where $h$ is the unknown spherical function to be found, $\nabla h$ and $\nabla^2 h$ denote the spherical gradient and spherical Hessian, and $I$ is the identity map. $\mathcal{D}$ is the duality operator, defined for $h \in C^{+}\left(S^{n-1}\right)$ by
$$
\mathcal{D} h(u)=\max _{v \in S^{n-1}} \frac{u \cdot v}{h(v)}, ~~ \text { for } u \in S^{n-1}.
$$

In Section \ref{SPre}, we introduce some basic concepts, such as convex bodies, radial Gauss  {map} and its reverse, $L_p$ dual curvature measures and $L_q$ dual intrinsic volume. In Section \ref{SLPCMADQ}, we define the $L_p$ curvature measures of ADQ and provided the variation formula on $\widetilde{C}_{p,m}^{\mathrm{a}}$. Moreover, in this section, we proved the continuity of the affine dual curvature measure $\widetilde{C}_{0,m}^{\mathrm{a}}$ and some basic properties of $\widetilde{C}_{p,m}^{\mathrm{a}}$.
In Section \ref{SSNSDm}, we solved the existence part of the discrete $L_p$ Minkowski problem on ADQ if $p > 1$ and proved the uniqueness of the solution for the case of polytopes when $p>mn$.
\begin{theorem}\label{T1.1}
 Let $m=1,\dots,n-1$,  $p>1$ and $p\neq mn$. Let $\mu$ be a discrete measure on $S^{n-1}$ that is not concentrated on any closed hemisphere. Then there exists a polytope $P\subset \mathbb{R}^n$ containing the origin in its interior such that $\widetilde{C}^a_{p,m}(P,\cdot)=\mu$.
\end{theorem}


In Section \ref{SSNSGm}, we turn to a general, possibly non-discrete Borel measure $\mu$ on $S^{n-1}$. We partially solve the existence part of the $L_p$ Minkowski problems.

\begin{theorem}\label{T1.3} Let $m=1,2,\dots,n-1$ and $p>n+m(n-2)$. Let $\mu$ be a finite Borel measure on $S^{n-1}$ that is not concentrated on any closed hemisphere. Then there exists a $K \in \mathcal{K}_{(o)}^n$  such that $ \widetilde{C}^a_{p,m}(K,\cdot)=\mu$.
\end{theorem}

 {In proving} Theorems \ref{T1.1} and \ref{T1.3}, we always need to estimate the value of the $m$ dimensional dual Radon transform $\mathcal{R}_m^*\left(|K \cap \cdot|^{n-1}\right)$.
To this end, we first estimate $\mathcal{R}_m^*\left(|B(x,r) \cap \cdot|^{n-1}\right)$  on the small ball $B(x,r)$ contained in $K$ and provide a specific expression for the estimate of $\mathcal{R}_m^*\left(|B(x,r) \cap \cdot|^{n-1}\right)$ using the radius of the small ball.

Finally in Section \ref{SSSm}, we consider the solution to the Minkowski problem for symmetric $L_p$ curvature measures of ADQ for $p\geq0$. We will demonstrate existence of solutions to the $L_p$ Minkowski problem on ADQ under the assumption that the given measure satisfies the strict subspace concentration inequality, i.e., for each proper subspace $\xi$ of $\mathbb{R}^n$,
\begin{equation}\label{E4.2}
\frac{\mu\left(\xi \cap S^{n-1}\right)}{|\mu|}<\frac{\operatorname{dim} \xi}{n}.
\end{equation}

\begin{theorem}\label{T1.2}
Suppose $p\geq0$, $m=1, \ldots, n-1$, and $\mu$ is a non-zero, even, finite Borel measure on $S^{n-1}$. If $\mu$ satisfies the strict subspace concentration inequality (\ref{E4.2}), then there exists an origin-symmetric convex body $K \in \mathcal{K}_{(o)}^n$ such that
$\widetilde{C}_{p,m}^{\mathrm{a}}(K, \cdot)=\mu$.
\end{theorem}

\section{Preliminaries}\label{SPre}
We develop some notations and, for quick later references, list some basic facts about convex
bodies. Good general references for the theory of convex bodies are provided by the books of Gardner \cite{Gardner}, Gruber \cite{Gruber} and Schneider \cite{Schneider14}.

We assume that $n\geq 2$ throughout the paper.  Let $\mathbb{R}^n$ denote $n$-dimensional Euclidean space with canonical inner product $x \cdot y$, for $x, y \in$ $\mathbb{R}^n$. Let $\mathbb{R}_{+}$ denote set of the positive real numbers. Let $o$ denote the origin of $\mathbb{R}^n.$
 Write $ {|x|}=\sqrt{x \cdot x}$ for the norm of $x\in\mathbb{R}^n$ and we normalize $x$ as $\langle x \rangle=x/ {| x |}$ for $x\in\mathbb{R}^{n}\backslash\{o\}$. For two vectors $x, y \in \mathbb{R}^n$, we denote by $[x, y]((x, y))$ the closed (open) segment with endpoints $x$ and $y$. Let $S^{n-1}$ denote the unit sphere in $\mathbb{R}^n$. Let $B^n(x_0,r_0)$ denote the ball centered at the $x_0$ with radius $r_0$ of $\mathbb{R}^n$ and $B^n$ denote the unit ball in $\mathbb{R}^n$ centered in origin. Let $C(S^{n-1})$ denote the set of continuous functions defined on $S^{n-1}$. Let $C_e(S^{n-1})$ and  $C^+(S^{n-1})$ denote the set of even, positive continuous functions, respectively. The group of special linear transformations of $\mathbb{R}^n$ is denoted by $SL(n)$. The $m$-dimensional Hausdroff measure is denoted by $\mathcal{H}^m$ and the $m$-dimensional Lebesgue measure is denoted by $\mathcal{L}^m$ . When considering an $m$-dimensional convex body $L$ that lies in an $m$-dimensional affine subspace of $\mathbb{R}^n$, we will write $|L|$ rather than $\mathcal{L}^m(L)$ to denote $m$-dimensional volume of $L$. Let $\xi^{\perp}$ denote the orthogonal complementary space of the linear subspace $\xi$ of $\mathbb{R}^n$. Let $u^{\perp}$ denote the orthogonal complementary space of the unit vector $u\in S^{n-1}$.  For a finite measure $\mu$ on $S^{n-1}$, we shall write $|\mu|$ for its total mass, that is $|\mu|=\mu\left(S^{n-1}\right)$.

\subsection{Convex bodies and star bodies}\label{s2.1}
 A {\it convex body} in $\mathbb{R}^n$ is a compact convex set. The set of convex bodies of $\mathbb{R}^n$ containing the origin is denoted by $\mathcal{K}_o^n$, and $\mathcal{K}^n_{(o)}$ denotes the set of convex bodies of $\mathbb{R}^n$  containing the origin in their interiors.  Let $\mathcal{K}_{(e)}^n$ denote the set of origin-symmetric convex bodies containing the origin in their interiors.  A {\it polytope} in $\mathbb{R}^n$ is the convex hull of a finite set of points in $\mathbb{R}^n$ provided it has positive volume
(i.e., $n$-dimensional volume). The convex hull of a subset of these points is called a {\it facet} of the polytope if it lies entirely on the boundary of the polytope and if it has positive area (i.e., $(n-1)$-dimensional volume).
Let $\operatorname{int}K$ denote the interior of $K\subset\mathbb{R}^n$ and $\partial K$ the boundary of $K$. For a convex body $K \in \mathcal{K}_{o}^n$ with $\operatorname{int}K\neq \emptyset$, we write $r(K)$ to denote the maximal radius of balls contained in $K$.

For two sets $K,L\subset\mathbb{R}^n$, the {\em Minkowski sum} $K+L$ is defined by
$K+L=\{x+y: x\in K,\;y\in L\}$.
The {\em scalar multiplication} $\lambda K$ of $K$ and $\lambda>0$ is defined as
$\lambda K=\{\lambda x:x\in K\}$.
The {\it support function}, $h(K,\cdot): S^{n-1} \rightarrow \mathbb{R}$, of $K \in \mathcal{K}^n$, is the continuous function on the unit sphere $S^{n-1}$, defined by
\begin{equation}\label{hKv}
h(K,v)=\max \{x \cdot v: x \in K\}. \nonumber
\end{equation}

A convex body is uniquely determined by its support function.
For $K\in\mathcal{K}^n$ and $u \in S^{n-1}$, the {\it support set} of $K$ with exterior unit normal $u$ is the set
$F(K, u)=\left\{x \in K: x\cdot u=h_K(u)\right\}$.
For $K\in\mathcal{K}_o^n$ with $\operatorname{int} K \neq \emptyset$, let
\begin{equation}\label{XiK}
{\color{red}\Xi_K=\left\{x \in \partial K:\; {\rm there\;exists\;exterior\;normal}\;\;v \in S^{n-1}\;{\rm at}\;x\;{\rm with}\;h_K(v)=0\right\}.}
\end{equation}

The {\em $L_p$ Minkowski sum} $K+_pL$ of $K,L\in \mathcal{K}_{o}^n$ and  $p> 0$ is defined by the support function
\begin{equation}
h(K+_pL,v)=\left(h(K,v)^p+h(L,v)^p\right)^{\frac{1}{p}}. \nonumber
\end{equation}
The collection of convex bodies can be viewed as a metric space with the {\it Hausdroff metric}:
\begin{eqnarray}
d_H(K,L):=\min\{t\geq 0:\;K\subset L+tB^n,\;L\subset K+tB^n\},\;\;K,L\in\mathcal{K}^n. \nonumber
\end{eqnarray}

The radial function, $\rho_K: x\in \mathbb{R}^n\backslash\{o\} \rightarrow \mathbb{R}$, of a compact star-shaped (with respect to the origin) set $K \subset \mathbb{R}^n$, is defined by
\begin{equation}\label{rhoK}
\rho_K(x)=\max \{t \geq 0: tx \in K\}.  \nonumber
\end{equation}
The star-shaped set $K$ is said to be a {\it star body}, if $\rho_K$ is continuous and positive. The collection of star bodies is denoted by $\mathcal{S}_{(o)}^n$.


\subsection{Radial Gauss image and its reverse.}\label{s2.2}

For $K\in\mathcal{K}_o^n$ and $z \in \partial K$, we write $N(K, z)$ to denote the {\it normal cone} at $z$; namely,
\begin{equation}
N(K, z)=\left\{y \in \mathbb{R}^n: y\cdot (x-z) \leq 0 \;\text { for } x \in K\right\}.\nonumber
\end{equation}

 Let $K \in \mathcal{K}_o^n$ with $\operatorname{int} K \neq \emptyset$. We recall that the so-called singular point $z \in \partial K$ where $\operatorname{dim} N(K, z) \geq 2$ form a Borel set of zero $\mathcal{H}^{n-1}$ measure, and hence its complement, the set of smooth points denoted by $\partial^{\prime} K$ is also a Borel set. For $z \in \partial^{\prime} K$, we write $\nu_K(z)$ to denote the unique exterior normal at $z$. In addition, for $z \in \partial K$, we define the {\it Gauss image} of $K$ at $z$ as
 $\boldsymbol{\nu}_K(z)=N(K, z) \cap S^{n-1}$.


Define the dual of $N(K, o)$ as
\begin{equation}\label{NKo}
N(K, o)^*=\left\{y \in \mathbb{R}^n: y\cdot x \leq 0 \text { for } x \in N(K, o)\right\}.
\end{equation}

For a $K \in \mathcal{K}_o^n$ with  $\operatorname{int} K \neq \emptyset$ and a Borel set $\eta \subset S^{n-1}$,   the  {\it reverse radial Gauss image} of $\eta$ is defined in \cite{BF19} as follows,
\begin{equation}
\boldsymbol{\alpha}_K^*(\eta)=\left\{u \in S^{n-1}: \rho_K(u) u \in \boldsymbol{\nu}_K^{-1}(\eta)\right\},\nonumber
\end{equation}
which extends the definition for $K\in\mathcal{K}^n_{(o)}$ in \cite{HLYZ16}.

For a $K \in \mathcal{K}_o^n$ with $\operatorname{int} K \neq \emptyset$ and $u \in S^{n-1}$,  the {\it radial map} of $K$ is defined by $r_K(u)=\rho_K(u) u \in \partial K$.
 For a Borel set $\eta \subset S^{n-1}$, the {\it radial Gauss image} of $K\in \mathcal{K}_o^n$ with $\operatorname{int}K\neq \emptyset$  is defined in \cite{BF19} as follows,
 \begin{equation}
 \boldsymbol{\alpha}_K(\eta)=\left\{ {v}\in S^{n-1}:\;{\rm there\;exists\;a}\; {u}\in\eta\;{\rm such\;that}\; {v}\in N(K, r_K( {u}))\right\},\nonumber
 \end{equation}
which extends the definition for $K\in\mathcal{K}^n_{(o)}$ in \cite{HLYZ16}.

For $\mathcal{H}^{n-1}$ almost all points $u \in S^{n-1} \cap\left(\operatorname{int} N(K, o)^*\right)$, there is a unique exterior unit normal $\alpha_K(u)$ at $r_K(u) \in \partial K$. Here $\alpha_K$ is the so-called {\it radial Gauss map}. For the other points $u \in S^{n-1} \cap\left(\operatorname{int} N(K, o)^*\right)$, we just choose an exterior unit normal $\alpha_K(u)$ at $r_K(u) \in \partial K$.

\subsection{Intersection
body}\label{s2.3}

The {\it intersection body} $\mathrm{I}K$ (first defined and named in \cite{Lutwak88}) of the star body $K \in \mathcal{S}_o^n$ can be defined by
\begin{equation}\label{DIB}
\rho_{\mathrm{I} K}(u)=\operatorname{vol}_{n-1}\left(K \cap u^{\perp}\right)=\frac{1}{n-1} \mathcal{R}\left(\rho_K^{n-1}\right)(u),
\end{equation}
where $\mathcal{R}: C\left(S^{n-1}\right) \rightarrow C\left(S^{n-1}\right)$ denotes the spherical Radon transform.

In order to study the generalized Busemann-Petty problem for $m$ dimensional sections, Zhang \cite{Zhang99} introduced $m$-intersection bodies. A star body $K$ is a {\it Zhang $m$-intersection body} if its radial function is the density of a measure obtained from the $m$-dimensional dual Radon transform $\mathcal{R}_m^*$. More precisely, $K$ is a Zhang $m$-intersection body if there is a measure $\mu$ on $G(n, m)$ such that $\mathcal{R}_m^* \mu$ is absolutely continuous with respect to spherical Lebesgue measure, and
$$
\rho_K(u)=\left(\mathcal{R}_m^* \mu\right)(d u) / d u.
$$
This definition extends the classical one since $\mathcal{R}: C\left(S^{n-1}\right) \rightarrow C\left(S^{n-1}\right)$ is self adjoint, and both $\mathcal{R}_{n-1}$ and $\mathcal{R}_{n-1}^*$ coincide with $\mathcal{R}$. See \cite{K00} for {\it Koldobsky generalized intersection bodies}, see also \cite{EM06,EM08}.

Recently, Cai-Leng-Wu-Xi \cite{CLWX} define the {\it bi-dual $m$-intersection body} $\mathbb{I}_mK$ of $K\in \mathcal{S}_o^n$ by
\begin{equation}\label{Dbib}
\rho_{\mathbb{I}_m K}(u)^{n-m}=\mathcal{R}_m^*\left(|K \cap \cdot|^{n-1}\right)(u) .
\end{equation}

Note that $\mathbb{I}_m K$ is always a Zhang $m$-intersection body. In fact, it can be defined by applying the $m$-dimensional Radon transform and its dual. In particular, $\mathbb{I}_{n-1}=c \mathbb{I}^2$, where $c=(n-1) /\left(n \omega_n\right)$ and $\mathbb{I}_1$ is radial symmetrization. Both the operators $\mathbb{I}_{n-1}$ and $\mathbb{I}_1$ are affine covariant and map $\mathcal{K}_{(e)}^n$ into $\mathcal{K}_{(e)}^n$.

\subsection{$L_p$ dual curvature measures and $L_q$ dual intrinsic volume}\label{s2.4}

The {\it $L_p$ dual curvature measures}, $\widetilde{C}_{p, q}$, are a two-parameter family of Borel measures on the unit sphere. Specifically, for $p, q \in \mathbb{R}$, a convex body $K \in \mathcal{K}_{(o)}^n$, and a star body $Q \in \mathcal{S}_{(o)}^n$, Lutwak-Yang-Zhang \cite{LYZ18} define the Borel measure $\widetilde{C}_{p, q}(K, Q, \cdot)$ on $S^{n-1}$ by
\begin{equation}
\int_{S^{n-1}} g(v) d \widetilde{C}_{p, q}(K, Q, v)=\frac{1}{n} \int_{S^{n-1}} g\left(\alpha_K(u)\right) h_K\left(\alpha_K(u)\right)^{-p} \rho_K(u)^q \rho_Q(u)^{n-q} d u,
\end{equation}
for each continuous $g: S^{n-1} \rightarrow \mathbb{R}$.

For simplicity of notation, we write $\widetilde{C}_{p, q}(K,u)$ instead of $\widetilde{C}_{p, q}(K, B^n,u)$.
According to Lemma 5.1 in Lutwak-Yang-Zhang \cite{LYZ18}, if $K \in \mathcal{K}_{(o)}^n$, $p \in \mathbb{R}$ and $q>0$, then for {all} Borel function $g: S^{n-1} \rightarrow \mathbb{R}$,
\begin{eqnarray}
\int_{S^{n-1}} g(u) d \widetilde{C}_{p, q}(K,u)&=&\frac{1}{n} \int_{S^{n-1}} g\left(\alpha_K(u)\right)\left(u\cdot\nu_K(\rho_K(u)u)\right)^{-p}\rho_K(u)^{q-p} du \label{E6.0}\\
&=&\frac{1}{n} \int_{\partial^{\prime} K} g\left(\nu_K(x)\right)(x\cdot \nu_K(x) )^{1-p} {|x|}^{q-n} d \mathcal{H}^{n-1}(x),\nonumber
\end{eqnarray}
{\color{red}where $\partial^{\prime}K$ denotes the set of points on $\partial K$ with unique unit outer normal vector.}

In \cite[Corollary 4.1]{BF19}, B\"or\"oczky-Fodor partially extend (\ref{E6.0}) to allow $o \in \partial K$:
 If $p>1, q>0$, $K \in \mathcal{K}_o^n$ with $\operatorname{int} K \neq \emptyset$, $\widetilde{C}_{p, q}\left(K, S^{n-1}\right)<\infty$ and $\mathcal{H}^{n-1}\left(\Xi_K\right)=0$ ( {see precise definition (\ref{XiK}) for $\Xi_K$}), then for every bounded Borel function $g: S^{n-1} \rightarrow \mathbb{R}$,
\begin{equation}\label{6.2}
\int_{S^{n-1}} g(u) d \widetilde{C}_{p, q}(K, u)=\frac{1}{n} \int_{\partial^{\prime} K} g\left(\nu_K(x)\right)(x\cdot \nu_K(x))^{1-p} {|x|}^{q-n} d \mathcal{H}^{n-1}(x).
\end{equation}

For $q>0$, B\"or\"oczky-Fodor \cite{BF19} extend Lutwak's definition of the $q$-th dual intrinsic volume
\begin{equation}\label{E6.11}
\widetilde{V}_q(K)=\frac{1}{n} \int_{S^{n-1}} \rho_K(u)^q du
\end{equation}
for $K \in \mathcal{K}_{(o)}^n$
to a compact convex set $K \in \mathcal{K}_o^n$ as
\begin{equation}\label{6.11}
\widetilde{V}_q(K)=\frac{1}{n} \int_{\partial^{\prime} K} (x\cdot \nu_K(x) ) {|x|}^{q-n} d \mathcal{H}^{n-1}(x).
\end{equation}

The following lemma provides a basic estimate on the inradius of $K$ in terms of its total $L_p$ dual curvature measure and $q$-th dual intrinsic
volume.
\begin{lemma}\cite[Lemma 6.10]{BF19}
For $n \geq 2$, $p>1$ and $q>0$, there exists a constant $c>0$ depending only on $p, q, n$ such that if $K \in \mathcal{K}_{(o)}^n$, then
\begin{equation}\label{25}
\widetilde{C}_{p, q}\left(K, S^{n-1}\right) \geq c \cdot r(K)^{-p} \cdot \widetilde{V}_q(K) .
\end{equation}
\end{lemma}

For $K\in\mathcal{K}_o^n$ with $\operatorname{int}K\neq \emptyset$, Lemma \ref{BF19} shows that $K$ contains the origin $o$ in its interior  when $K$ can be approximated by a sequence of convex bodies $K_j\in\mathcal{K}^n_{(o)}$.

%

\begin{lemma}\cite[Lemma 6.11]{BF19}\label{BF19} {Suppose $p>1$ and $0<q \leq p$. If a sequence of convex bodies $K_j \in \mathcal{K}_{(o)}^n$ converges to a convex body $K \in \mathcal{K}_o^n$ with non-empty interior, and if the sequence of dual curvature measures $\widetilde{C}_{p, q}(K_j, S^{n-1})$ remains bounded, then it follows that $K$ must also belong to $\mathcal{K}_{(o)}^n$.}
\end{lemma}


\section{$L_p$ curvature measures of ADQ on $\mathcal{K}^n_{o}$}\label{SLPCMADQ}
\subsection{Definitions of $L_p$ curvature measures of ADQ}
In this section, we extend the definition of the affine dual curvature measure $\widetilde{C}^a_m(K, \cdot)$ for $K \in \mathcal{K}_{(o)}^n$ in \cite{CLWX}  to  the $L_p$ curvature measure of ADQ denoted by $ \widetilde{C}^a_{p,m}(K, \cdot)$ for $K \in \mathcal{K}_o^n$.

In \cite{CLWX}, the {\it affine dual curvature measure} of $K\in\mathcal{K}_{(o)}^n$ is defined  by
\begin{equation}\label{2.4}
\widetilde{C}_{m}^a(K, \eta)=\int_{\boldsymbol{\alpha}_K^*(\eta)}\rho^{m}_K(u)\mathcal{R}_m^*\left(\left|K \cap \cdot\right|^{n-1}\right)(u)d u,
\end{equation}
for each Borel set $\eta \subset S^{n-1}$.

  Let $p\in\mathbb{R}$, $m=1, \ldots, n-1$, and  $K \in \mathcal{K}_{o}^n$, we define the {\it $L_p$ curvature measure of ADQ} $\widetilde{C}_{p,m}^a(K, \cdot)$ of $K$ for each Borel set $\eta \subset S^{n-1}$ by,
\begin{eqnarray}
&&\widetilde{C}^a_{p,m}(K, \eta)=0\;\;\;{\rm if}\;\;\operatorname{dim} K \leq n-1,\label{E22}\\
&&\widetilde{C}_{p,m}^a(K, \eta)=\int_{ {\boldsymbol{\alpha}_K^*(\eta)}} \left(u\cdot\nu_K\left(\rho_K(u)u\right)\right)^{-p}\rho_K(u)^{m-p} \mathcal{R}_m^*\left(\left|K \cap \cdot\right|^{n-1}\right)(u) d u \label{E2.2}
\end{eqnarray}
{\color{red}if $\operatorname{int} K \neq \emptyset$  and $\int_{S^{n-1}}\left(u\rho_K(u)\cdot\nu_K\left(\rho_K(u)u\right)\right)^{-p}du<\infty$.}

If $K \in \mathcal{K}_{(o)}^n$, then for  each Borel set $\eta \subseteq S^{n-1}$,
\begin{equation}\label{2.3}
\widetilde{C}_{p,m}^a(K, \eta)=\int_{\boldsymbol{\alpha}_K^*(\eta)} \left(u\cdot\nu_K\left(\rho_K(u)u\right)\right)^{-p}\rho_K(u)^{m-p} \mathcal{R}_m^*\left(\left|K \cap \cdot\right|^{n-1}\right)(u)d u,
\end{equation}
or equivalently,
\begin{equation}\label{E3.12}
\widetilde{C}^a_{p,m}(K,\eta)=\int_{\boldsymbol{\nu}^{-1}_K(\eta)}\left(x \cdot \nu_K(x)\right)^{1-p} {|x|}^{m-n}
\mathcal{R}_m^*\left(\left|K \cap \cdot\right|^{n-1}\right)(\langle x\rangle)d \mathcal{H}^{n-1}(x).
\end{equation}
 {See Section \ref{SPre} for the definition of $\langle x\rangle.$ }%

By (\ref{E3.12}) and (\ref{E6.0}), for $K \in \mathcal{K}_{(o)}^n$, we have
\begin{equation}\label{E25}
d\widetilde{C}_{p,m}^a(K, \cdot)=h_K^{-p}(\cdot)d \widetilde{C}^a_{m}(K, \cdot).
\end{equation}

\subsection{The variation formula on $\widetilde{C}_{p,m}^a$}

 Let $h_0\in C^+(S^{n-1})$ and $f\in C(S^{n-1})$. For sufficiently small $\delta>0$ and $t \in(-\delta, \delta)$, let
\begin{eqnarray}\label{3.1}
h_t(v)= \begin{cases}\left(h^p_0(v)+tf(v)+o(t,v)\right)^{\frac{1}{p}}, & \text { if } p >0, \\ h_0(v)e^{tf(v)+o(t,v)}, & \text { if } p=0,\end{cases}
\end{eqnarray}
where $\lim _{t \rightarrow 0} o(t, \cdot) / t=0$ uniformly on $S^{n-1}$.
The {\it $L_p$ family of Wulff shapes} $K_t$ associated with $\left(h_0, f\right)$ is defined by
\begin{equation}\label{3.2}
K_t=\left\{x \in \mathbb{R}^n: x \cdot v \leq h_t(v) \text { for all } v \in S^{n-1}\right\}, \quad t \in(-\delta, \delta).
\end{equation}

{\color{red}The proof of the following lemma is similar to that of \cite[Theorem 3.1]{CLWX}. It is important to note, however, that the result in \cite{CLWX} does not cover the general case involving the term
$o(t,v)$ that we consider here.}


\begin{lemma}\label{CLWX} Let $p=0$ and $m=1,\dots,n-1$, and let $K_t$ be the  {log-family} of Wulff shapes defined by (\ref{3.2}). Then,
\begin{equation}
\lim _{t \rightarrow 0} \frac{\widetilde{\Psi}_m\left(K_t\right)-\widetilde{\Psi}_m(K_0)}{t}=n \int_{S^{n-1}} f(u) d \widetilde{C}_m^{\mathrm{a}}(K_0, u).\nonumber
\end{equation}
\end{lemma}

Using Lemma \ref{CLWX} and the definition of the $L_p$ family of Wulff shapes, $L_p$ curvature measure of ADQ can be derived from the following variation formula on $\widetilde{\Psi}_m\left(K_t\right)$ at $t=0$.

\begin{theorem}\label{2.1}
Let $p>0$ and $m=1,\dots,n-1$, and let $K_t$ be the $L_p$ family of Wulff shapes defined by (\ref{3.2}). Then,
\begin{equation}\label{deq}
\lim _{t \rightarrow 0} \frac{\widetilde{\Psi}_m\left(K_t\right)-\widetilde{\Psi}_m(K_0)}{t}=\frac{n}{p} \int_{S^{n-1}} f(u) d \widetilde{C}_{p,m}^a(K_0, u).
\end{equation}
\end{theorem}

\begin{proof}
When $p>0$, by (\ref{3.1}),
\begin{eqnarray}
h_t(u)&=&h_0(u)\left(1+t\frac{f(u)}{h_0(u)^p}+o(t,u)\right)^{\frac{1}{p}}.\nonumber
\end{eqnarray}
Thus,
\begin{eqnarray}\label{loghpt}
\log h_t(u)&=&\log h_o(u)+\frac{1}{p}\log \left(1+t\frac{f(u)}{h_o(u)^p}+o(t,u)\right)\nonumber\\
&=&\log h_o(u)+\frac{t}{p}\frac{f(u)}{h_o(u)^p}+o(t,u).
\end{eqnarray}
By (\ref{loghpt}) and Lemma \ref{CLWX},
\begin{equation}\label{Elim}
\lim _{t \rightarrow 0} \frac{\widetilde{\Psi}_m\left(K_t\right)-\widetilde{\Psi}_m(K_0)}{t}=\frac{n}{p} \int_{S^{n-1}} \frac{f(u)}{h_0(u)^p} d \widetilde{C}_{m}^a(K_0, u).
\end{equation}
Since $h_0(u)=h_{K_0}(u)$ for almost every $u\in S^{n-1}$ on surface measure $dS(K_0,\cdot)$, and by (\ref{E25})  and (\ref{Elim}), we can get the desired equation (\ref{deq}).
\end{proof}

The following corollary  shows the relation between $L_p$ curvature measure of ADQ   and the variation formula of $\widetilde{\Psi}_m\left(K+_pt\cdot L\right)$ on  $K,L\in\mathcal{K}^n_{(o)}$ at $t=0$.

\begin{corollary}\label{lpaffinedualmeasuare}
Let $p>0$ and $m=1,\dots,n-1$, and let $K,L\in\mathcal{K}^n_{(o)}$. Then,
\begin{equation}\label{difflpaffinedualmeasuare}
\lim _{t \rightarrow 0} \frac{\widetilde{\Psi}_m\left(K+_pt\cdot L\right)-\widetilde{\Psi}_m(K)}{t}=\frac{n}{p} \int_{S^{n-1}}h_L(u)^p d \widetilde{C}_{p,m}^a(K, u).
\end{equation}
\end{corollary}

The following corollary is the discrete version of Theorem \ref{2.1}.

\begin{corollary}\label{C2.3}
 Let $p>0$ and $m=1,\dots,n-1$, and  {$u_1, \ldots, u_N$} be $N$ unit vectors that are not contained in any closed hemisphere and $z=\left(z_1, \ldots, z_N\right) \in\left(\mathbb{R}_{+}\right)^N$. Let $\beta=\left(\beta_1, \ldots, \beta_N\right) \in \mathbb{R}^N$. For
 sufficiently small $|t|$ and a higher-order infinitesimal $o(t)$, define $z(t)$ as follows,
\begin{equation}\label{3.4}
\left(z(t)\right)_i=\left(z_i^p+t\beta_i+o(t)\right)^{\frac{1}{p}},\; i=1,\dots,N,
\end{equation}
and
\begin{equation}
P_t=[z(t)]=\bigcap_{i=1}^N\left\{x \in \mathbb{R}^n: x \cdot {u_i} \leq (z(t))_i\right\}.\nonumber
\end{equation}
Then,
\begin{equation}
\lim _{t \rightarrow 0} \frac{\widetilde{\Psi}_m\left(P_t\right)-\widetilde{\Psi}_m(P_0)}{t}=\frac{n}{p}\sum_{i=1}^N \beta_i\widetilde{C}_{p,m}^a(P_0, \{ {u_i}\}).\nonumber
\end{equation}
\end{corollary}

\subsection{The continuity of affine dual curvature measures}

In this section, we prove that the affine dual quermassintegral $\widetilde{\Psi}_m$ is a continuous function of $K \in \mathcal{K}_o^n$  {in the Hausdorff metric}. Moreover, using the continuity of $\widetilde{\Psi}_m$, we proved that the affine dual curvature measure $\widetilde{C}^a_{m}\left(K,\cdot\right)$ is weakly convergent on $K\in\mathcal{K}_o^n$.

%

\begin{lemma}\label{LPsi}
If  $K\in\mathcal{K}^n_{o}$, then
\begin{equation}\label{E4.4}
\widetilde{\Psi}_m(K)=\frac{1}{m}\int_{S^{n-1}}\rho_K(u)^m\mathcal{R}_m^*\left(\left|K \cap \cdot\right|^{n-1}\right)(u) du.
\end{equation}
\end{lemma}
\begin{proof}
For $K\in\mathcal{K}^n_{(o)}$, let $f\equiv 1$ and $h_0=h_K$ in Lemma \ref{CLWX} {. Then,} we have
\begin{equation}
n \int_{S^{n-1}} 1 d \widetilde{C}_m^{\mathrm{a}}(K, u)=\lim _{t \rightarrow 0} \frac{\widetilde{\Psi}_m\left(e^tK\right)-\widetilde{\Psi}_m(K)}{t}=mn\widetilde{\Psi}_m(K).\nonumber
\end{equation}
Therefore, by (\ref{2.4}),
\begin{equation}\label{EPsimK}
\widetilde{\Psi}_m(K)=\frac{1}{m}\int_{S^{n-1}}1 d \widetilde{C}_m^{\mathrm{a}}(K, u)=\frac{1}{m}\int_{S^{n-1}}\rho_K(u)^{m} \mathcal{R}_m^*\left(\left|K \cap \cdot\right|^{n-1}\right)(u) d u.
\end{equation}

If $K\in\mathcal{K}^n_{o}$, then there exits a sequence  $\{K_j\}_{j\in\mathbb{N}}\subset \mathcal{K}^n_{(o)}$ tending to $K$ with respect to Hausdorff distance. By the definition  (\ref{2.2}) of $\widetilde{\Psi}_m(K)$, we have
\begin{eqnarray}\label{limPsim}
\widetilde{\Psi}_m(K)=\int_{G(n, m)} \left|K \cap \xi\right|^n d \nu_m(\xi)
&=&\frac{1}{m^n}\int_{G(n, m)}\left(\int_{S^{n-1}\cap\xi}\rho_K(u)^mdu\right)^n d \nu_m(\xi)\\
&=&\lim_{j\rightarrow\infty}\frac{1}{m^n}\int_{G(n, m)}\left(\int_{S^{n-1}\cap\xi}\rho_{K_j}(u)^mdu\right)^n d \nu_m(\xi)\nonumber\\
&=&\lim_{j\rightarrow\infty}\widetilde{\Psi}_m(K_j).\nonumber
\end{eqnarray}

By (\ref{limPsim}) and (\ref{EPsimK}),
\begin{eqnarray}
\widetilde{\Psi}_m(K)=\lim_{j\rightarrow\infty}\widetilde{\Psi}_m(K_j)
&=&\lim_{j\rightarrow\infty}\frac{1}{m}\int_{S^{n-1}}\rho_{K_j}(u)^{m} \mathcal{R}_m^*\left(\left|K_j \cap \cdot\right|^{n-1}\right)(u) d u\nonumber\\
&=&\frac{1}{m}\int_{S^{n-1}}\rho_{K}(u)^{m} \mathcal{R}_m^*\left(\left|K \cap \cdot\right|^{n-1}\right)(u)  d u.\nonumber
\end{eqnarray}
This is the desired equation (\ref{E4.4}).
\end{proof}

 {Note that we only prove the formula for $\widetilde{\Psi}_m$ in $\mathcal{K}_o^n$. Next, we prove the continuity of $\widetilde{\Psi}_m$ in $\mathcal{K}_o^n$.
}

\begin{lemma}\label{L4.1}
 For $m=1,2,\dots,n-1$, $\widetilde{\Psi}_m(K)$ is a continuous function of $K \in \mathcal{K}_o^n$  {in the Hausdorff metric}.
\end{lemma}
\begin{proof}
Let $R>0$ be such that $K \subset \operatorname{int} R B^n$. Let the sequence of $\{K_j\}_{j\in\mathbb{N}} \subset \mathcal{K}_o^n$ converge to $K$ {\color{red}in Hausdorff metric.} Thus, we can assume that $K_j \subset R B^n$ for all $K_j$.

If $\operatorname{dim} K \leq n-1$, then there exists $v \in S^{n-1}$ such that $K \subset v^{\perp}$. In this case, by (\ref{E4.4}), $\widetilde{\Psi}_m(K)=0$.
For $t \in[0,1)$, let
$$
\Xi(v, t)=\left\{x \in \mathbb{R}^n:|v\cdot x| \leq t\right\}.
$$
There exists a $t_0 \in(0,1)$ such that for  {all} $t \in\left(0, t_0\right)$ and $v \in S^{n-1}$ it holds that
\begin{equation}\label{E4.8}
\mathcal{H}^{n-1}\left(S^{n-1} \cap \Xi(v, t)\right)<3 t(n-1) \omega_{n-1}.
\end{equation}
Let $\varepsilon \in\left(0, t_0\right)$. Since $K_j \rightarrow K$ in the Hausdorff metric, there exists a sufficiently large $j_{\varepsilon}\in\mathbb{N}$ such that for all $j>j_{\varepsilon}$, it holds that
$$K_j \subset K+\varepsilon^2 B^n \subset \Xi\left(v, \varepsilon^2\right).$$
 Then if $u \in S^{n-1} \backslash \Xi(v, \varepsilon)$, we have
$$
\varepsilon^2 \geq\left(v\cdot\rho_{K_j}(u) u\right)=\rho_{K_j}(u) (v\cdot u) \geq \rho_{K_j}(u) \cdot \varepsilon.
$$
Thus, we have
 \begin{equation}\label{31}
\rho_{K_j}(u) \leq \varepsilon,\; {\rm for\; {all}}\;u \in S^{n-1} \backslash \Xi(v, \varepsilon).
\end{equation}

We deduce from (\ref{E4.4}), (\ref{E4.8}), (\ref{31}), (\ref{Rmstar}) and $K_j \subset R B^n$ that for  {all} $\varepsilon \in\left(0, t_0\right)$, if $j>j_{\varepsilon}$, then
$$
\begin{aligned}
\widetilde{\Psi}_m\left(K_j\right) & \leq \frac{1}{m}\int_{S^{n-1} \backslash \Xi(v, \varepsilon)} \varepsilon^m \mathcal{R}_m^*\left(\left|K_j \cap \cdot\right|^{n-1}\right)(u) d \mathcal{H}^{n-1}(u)\\
&+\frac{1}{m}\int_{S^{n-1} \cap \Xi(v, \varepsilon)} R^m   \mathcal{R}_m^*\left(\left|K_j \cap \cdot\right|^{n-1}\right)(u)  d \mathcal{H}^{n-1}(u) \\
& \leq \omega_m^nR^{m(n-1)} \varepsilon^m+\frac{3(n-1) \omega_m^n\omega_{n-1} R^{mn}}{n\omega_n}\varepsilon.
\end{aligned}
$$
Therefore $\lim _{j \rightarrow \infty} \widetilde{\Psi}_m\left(K_j\right)=0=\widetilde{\Psi}_m(K)$.

Next, we consider the case of $\operatorname{int}K \neq \emptyset$ and $o\in\partial K$.
By (\ref{E4.4}), $\mathcal{H}^{n-1}\left(S^{n-1} \cap \partial N(K, o)^*\right)=0$ and Lebesgue's dominated convergence theorem,
\begin{align*}
\lim_{j\rightarrow\infty} \widetilde{\Psi}_m(K_j)= &\lim_{j\rightarrow\infty} \frac{1}{m}\int_{S^{n-1}}\rho_{K_j}^m(u) \mathcal{R}_m^*\left(\left|K_j \cap \cdot\right|^{n-1}\right)(u)  du \\
  = &\lim_{j\rightarrow\infty} \frac{1}{m}\int_{S^{n-1}\backslash \partial N(K_j, o)^*}\rho_{K_j}^m(u) \mathcal{R}_m^*\left(\left|K_j \cap \cdot\right|^{n-1}\right)(u) du \\
  = &\frac{1}{m}\int_{S^{n-1}\backslash \partial N(K, o)^*}\rho_{K}^m(u) \mathcal{R}_m^*\left(\left|K \cap \cdot\right|^{n-1}\right)(u) du=\widetilde{\Psi}_m(K).
\end{align*}

Finally, let $K \in\mathcal{K}^n_{(o)}$. By \cite[Theorem 2.2]{LW23}, if
$d_H(K_j,K)\rightarrow 0$, then $\left|\rho_{K_j}(u)-\rho_K(u)\right|\rightarrow0$ uniformly when $j\rightarrow \infty$.
 Thus, by (\ref{E4.4}), we have $\lim _{j \rightarrow \infty} \widetilde{\Psi}_m\left(K_j\right)=\widetilde{\Psi}_m(K)$.
\end{proof}

\begin{proposition}\label{P2.6}
 If $m=1,\dots,n-1$, and $\left\{K_j\right\}, j \in \mathbb{N}$, converges to $K$ for $K_j, K \in \mathcal{K}_o^n$, then $\widetilde{C}^a_{m}\left(K_j,\cdot\right)$ converges weakly to $\widetilde{C}^a_{m}(K, \cdot)$.
 \end{proposition}

 \begin{proof}
Since any element of $\mathcal{K}_o^n$ can be approximated by elements of $\mathcal{K}_{(o)}^n$, we may assume that each $K_j \in \mathcal{K}_{(o)}^n$. We fix $R>0$ such that $K \subset \operatorname{int} R B^n$, and hence we may also assume that $K_j \subset R B^n$ for all $K_j$. We need to prove that if $g: S^{n-1} \rightarrow \mathbb{R}$ is continuous, then
\begin{equation}\label{E33}
\lim _{j \rightarrow \infty} \int_{S^{n-1}} g(u) d \widetilde{C}^a_{m}\left(K_j, u\right)=\int_{S^{n-1}} g(u) d \widetilde{C}^a_{m}(K, u).
\end{equation}

If $\operatorname{dim} K \leq n-1$, then $\widetilde{C}^a_m(K, \cdot)$ is the constant zero measure by (\ref{E22}). By $\widetilde{C}^a_m\left(K_j, S^{n-1}\right)=\widetilde{\Psi}_m\left(K_j\right)$, $\widetilde{\Psi}_m\left(K\right)=0$ and Lemma \ref{L4.1},
\begin{eqnarray}
\lim _{j \rightarrow \infty} \int_{S^{n-1}} g(u) d \widetilde{C}^a_m\left(K_j, u\right)&\leq& \lim _{j \rightarrow \infty} \int_{S^{n-1}} \max\{g(u):\;u\in S^{n-1}\} d \widetilde{C}^a_m\left(K_j, u\right)\nonumber\\
&=& \max\{g(u):\;u\in S^{n-1}\}\lim _{j \rightarrow \infty}  \widetilde{C}^a_m\left(K_j, S^{n-1}\right)\nonumber\\
&=&\max\{g(u):\;u\in S^{n-1}\}\lim _{j \rightarrow \infty} \widetilde{\Psi}_m\left(K_j\right)\nonumber\\
&=&\max\{g(u):\;u\in S^{n-1}\}\widetilde{\Psi}_m\left(K\right)\nonumber\\
&=&0=\int_{S^{n-1}} g(u) d \widetilde{C}^a_m(K, u).\nonumber
\end{eqnarray}

Next, we consider the case  of $\operatorname{int} K \neq \emptyset$ and $o\in\partial K$. By (\ref{2.4}) and (\ref{E2.2}), (\ref{E33}) is equivalent to
\begin{eqnarray}\label{E34}
&&\lim _{j \rightarrow \infty} \int_{S^{n-1}} g\left(\alpha_{K_j}(u)\right) \rho_{K_j}(u)^m \mathcal{R}_m^*\left(\left|K_j \cap \cdot\right|^{n-1}\right)(u) d \mathcal{H}^{n-1}(u)\\
&=&\int_{S^{n-1} \cap(\mathrm{int}N(K, o)^*)} g\left(\alpha_K(u)\right) \rho_K(u)^m \mathcal{R}_m^*\left(\left|K \cap \cdot\right|^{n-1}\right)(u)d \mathcal{H}^{n-1}(u).\nonumber
\end{eqnarray}

Since $\langle \cdot\rangle_K$ is Lipschitz and $\mathcal{H}^{n-1}\left(S^{n-1} \cap(\partial N(K,o)^{\ast})\right)=0$, to prove (\ref{E34}),   it suffices to show
\begin{eqnarray}\label{E35}
&&\lim _{j \rightarrow \infty} \int_{\langle \cdot\rangle_K\left((\operatorname{int} N(K, o)^*) \cap \partial^{\prime} K\right)} g\left(\alpha_{K_j}(u)\right) \rho_{K_j}(u)^m \mathcal{R}_m^*\left(\left|K_j \cap \cdot\right|^{n-1}\right)(u)d \mathcal{H}^{n-1}(u)\\
 &=&\int_{\langle \cdot\rangle_K\left((\operatorname{int} N(K, o)^*) \cap \partial^{\prime} K\right)} g\left(\alpha_K(u)\right) \rho_K(u)^m \mathcal{R}_m^*\left(\left|K \cap \cdot\right|^{n-1}\right)(u)  d \mathcal{H}^{n-1}(u)\nonumber
\end{eqnarray}
and
\begin{equation}\label{E36}
\lim _{j \rightarrow \infty} \int_{S^{n-1} \backslash N(K, o)^*} g\left(\alpha_{K_j}(u)\right) \rho_{K_j}(u)^m\mathcal{R}_m^*\left(\left|K \cap \cdot\right|^{n-1}\right)(u) d \mathcal{H}^{n-1}(u)=0.
\end{equation}

For (\ref{E35}), let $u \in \langle \cdot\rangle_K\left((\operatorname{int} N(K, o)^*) \cap \partial^{\prime} K\right)$ ( {See Section \ref{SPre} for the definition of $\langle \cdot\rangle_K.$}). Since that $K_j$ converges to $K$ on Hausdorff metric, we have
\begin{equation}\label{varU}
\lim _{j \rightarrow \infty} \rho_{K_j}(u)^m=\rho_K(u)^m\;{\rm and}\;\lim_{j\rightarrow\infty}\mathcal{R}_m^*\left(\left|K_j \cap \cdot\right|^{n-1}\right)(u)=\mathcal{R}_m^*\left(\left|K \cap \cdot\right|^{n-1}\right)(u).
\end{equation} Since $\alpha_K(u)$ is the unique normal at $\rho_K(u) u \in \partial^{\prime} K$, we have $\lim _{j \rightarrow \infty} \alpha_{K_j}(u)=\alpha_K(u)$. Moreover, by the continuity of $g$, we have
\begin{equation}\label{alphalim}
\lim _{j \rightarrow \infty} g\left(\alpha_{K_j}(u)\right)=g\left(\alpha_K(u)\right).
\end{equation}
By Lebesgue's dominated convergence theorem, (\ref{E35})  follows from (\ref{varU}), (\ref{alphalim}) and $K_j\subset RB^n$.

To prove (\ref{E36}), by Lebesgue's dominated convergence theorem,  we only need to show that
\begin{equation}\label{5.19}
\lim_{j\rightarrow\infty}\rho_{K_j}(u)=0,\;\; u\in S^{n-1} \backslash N(K, o)^*.
\end{equation}
Since $u \notin N(K, o)^*$, there exists $v \in N(K, o)$ such that $ v\cdot u=\delta>0$. Let $\varepsilon>0$. As $h_K(v)=0$ and $K_j$ converges to $K$, there exists $j_0$ such that $h_{K_j}(v) \leq \delta \varepsilon$ if $j \geq j_0$. In particular, if $j \geq j_0$, then
$$
\varepsilon \delta \geq h_{K_j}(v) \geq \rho_{K_j}(u)(v\cdot u)=\rho_{K_j}(u) \delta,
$$
which implies that (\ref{5.19}) is established, and consequently (\ref{E36}) is  established.

Finally, when $o \in \operatorname{int} K$, we can get the desired conclusion by the same proof as  (\ref{E35}).
 \end{proof}
\subsection{Some properties on $L_p$ curvature measures of ADQ}
In this section, we will prove a basic estimate on the inradius of $K\in\mathcal{K}^n_{(o)}$ in terms of its total $L_p$ curvature measure of ADQ (see Lemma \ref{L6.2}). Here the outradius of $K\in\mathcal{K}^n_{o}$, denoted by $R\left(K\right)$, denotes the  minimal radius of balls centered at $o$ and containing $K$.
 Moreover, we will prove that $K$ contains the origin $o$ in its interior if $K\in\mathcal{K}^n_{o}$  with $\operatorname{int}K \neq \emptyset$  and  $K$ can be approximated by a sequence of convex bodies $K_j\in\mathcal{K}^n_{(o)}$ (see Lemma \ref{L6.4}).

\begin{lemma}\label{L5.6}
Let $n>m\geq 2$, $x_0\in\mathbb{R}^n$, $r_0,R_0\in\mathbb{R}_+$ satisfy $B^n(x_0,r_0)\subset B^n(o,R_0)$. Then there exists a constant $c_2>0$ depending on $n,m,R_0$ such that for  {all} $u\in S^{n-1}$,
\begin{equation}\label{c0}
\mathcal{R}_m^*\left(\left|B(x_0,r_0) \cap \cdot\right|^{n-1}\right)(u)\geq c_2 r_0^{n+m(n-2)}.
\end{equation}
\end{lemma}
\begin{proof}
Translating $B^n(x_0,r_0)$ along the straight line parallel to $u$ and passing through point $x_0$ to $B^n(x_0+t_0u,r_0)$ such that $x_0+t_0u\in u^{\perp}$. Let $v_0=\left(x_0+t_0u\right)/\|x_0+t_0u\|$ if $x_0+t_0u\neq o$. Otherwise, choose arbitrarily $v_0\in S^{n-1}\cap u
^\perp$. Translating $B^n(x_0+t_0u,r_0)$ along the straight line parallel to $v_0$ and passing through point $o$ to
$B^n(x_0+t_0u+t_1v_0,r_0)$ such that
$\|x_0+t_0u+t_1v_0\|=R_0-r_0$. Thus,
\begin{eqnarray}\label{Gnmo}
\mathcal{R}_m^*\left(\left|B(x_0,r_0) \cap \cdot\right|^{n-1}\right)(u)
&=&\mathcal{R}_m^*\left(\left|B(x_0+t_0u,r_0) \cap \cdot\right|^{n-1}\right)(u)\\
&\geq& \mathcal{R}_m^*\left(\left|B(x_0+t_0u+t_1v_0,r_0) \cap \cdot\right|^{n-1}\right)(u).\nonumber
\end{eqnarray}

By the definition (\ref{Rmstar}) of $\mathcal{R}_m^*$, there exits a constant $c_2>0$  depending on $n,m,R_0$ such that
\begin{equation}\label{t0t1}
\mathcal{R}_m^*\left(\left|B(x_0+t_0u+t_1v_0,r_0) \cap \cdot\right|^{n-1}\right)(u)\geq c_2 r_0^{n+m(n-2)}.
\end{equation}

Combining (\ref{Gnmo}) with (\ref{t0t1}), we can get the desired conclusion.
\end{proof}

%

\begin{lemma}\label{L6.2}
 For $p>1$ and $m=1,\dots,n-1$, if $K \in \mathcal{K}_{(o)}^n$ and $R(K)<R_0$, then there exists a constant $c_0>0$  depending $n,m,p$ and $R_0$ such that
\begin{equation}\label{E49}
\widetilde{C}^a_{p, m}\left(K, S^{n-1}\right) \geq c_0r(K)^{-p+n+m(n-2)}   \widetilde{\Psi}_m(K).
\end{equation}
\end{lemma}
\begin{proof}
  By (\ref{E4.4}), (\ref{E6.11}),  $R(K)<R_0$ and the definition (\ref{Rmstar}) of $\mathcal{R}_m^*$, we have
  \begin{eqnarray}\label{Psi1}
  \widetilde{\Psi}_m(K)&=&\frac{1}{m}\int_{S^{n-1}}\rho_K(u)^m\mathcal{R}_m^*\left(\left|K \cap \cdot\right|^{n-1}\right)(u) du\\
  &\leq&\frac{\omega^n_m}{\omega_n}R_0^{m(n-1)}\widetilde{V}_m(K).\nonumber
  \end{eqnarray}

  By (\ref{2.3}), (\ref{E6.0}) and (\ref{c0}), there exists a constant $c_0>0$ depending on $n,m,R_0$ such that
  \begin{eqnarray}\label{Psi2}
 \widetilde{C}^a_{p, m}\left(K, S^{n-1}\right)&=& \int_{S^{n-1}} \left(u\cdot\nu_K\left(\rho_K(u)u\right)\right)^{-p}\rho_K(u)^{m-p} \mathcal{R}_m^*\left(\left|K \cap \cdot\right|^{n-1}\right)(u)d u\nonumber\\
 &\geq&c_0r_0^{n+m(n-2)}\widetilde{C}_{p, m}\left(K, S^{n-1}\right).
  \end{eqnarray}

  Combining (\ref{Psi1}), (\ref{Psi2}) and (\ref{25}), there exists a constant $c_0>0$ depending on $n,m,R_0$ such that
  \begin{equation*}
\widetilde{C}^a_{p, m}\left(K, S^{n-1}\right) \geq c_0r_0^{-p+n+m(n-2)}\widetilde{\Psi}_m(K),
\end{equation*}
which is the desired conclusion.
\end{proof}

\begin{lemma}\label{L6.4}
If $p>1$, $0<q \leq p$ and $K_j \in \mathcal{K}_{(o)}^n$ for $j \in \mathbb{N}$ converge to $K \in \mathcal{K}_o^n$ with $\operatorname{int} K \neq \emptyset$ such that $\widetilde{C}^a_{p,m}\left(K_j, S^{n-1}\right)$ stays bounded, then $K \in \mathcal{K}_{(o)}^n$.
\end{lemma}
\begin{proof}
Since $K_j \in \mathcal{K}_{(o)}^n$ for $j \in \mathbb{N}$ converge to $K \in \mathcal{K}_o^n$ and $\operatorname{int} K \neq \emptyset$, there exists a ball $B(x_0,r_0)\subset K$ and $B(x_0,r_0)\subset K_j$ and there exists a ball $B(o,R_0)$ such that $K\subset B(o,R_0)$ and $K_j\subset B(o,R_0)$, $j\in \mathbb{N}$.

 If $n>m\geq 2$, by Lemma \ref{L5.6}, there exists a constant $c_2$ depending on $n,m,R_0$ such that for  {all} $u\in S^{n-1}$,
\begin{equation}\label{geqc0}
\mathcal{R}_m^*\left(\left|B(x_0,r_0) \cap \cdot\right|^{n-1}\right)(u)\geq c_2 r_0^{n+m(n-2)}.
\end{equation}
%
%

By (\ref{6.2}), (\ref{E3.12}) and (\ref{geqc0}), we have
\begin{equation}
\widetilde{C}^a_{p,m}\left(K_j, S^{n-1}\right)\geq nc_2 r_0^{n+m(n-2)}\widetilde{C}_{p,m}\left(K_j, S^{n-1}\right),\nonumber
\end{equation}
which implies that $\widetilde{C}_{p,m}\left(K_j, S^{n-1}\right)$ stays bounded if $\widetilde{C}^a_{p,m}\left(K_j, S^{n-1}\right)$ stays bounded. From Lemma \ref{BF19}, we obtain $K \in \mathcal{K}_{(o)}^n$.
%

If $n\geq 2$ and $m=1$, then by (\ref{2.3}) and (\ref{E6.0}), we have
\begin{eqnarray}
\widetilde{C}^a_{p,1}\left(K_j, S^{n-1}\right)&=&\int_{S^{n-1}} \left(u\cdot\nu_{K_j}\left(\rho_{K_j}(u)u\right)\right)^{-p}\rho_{K_j}(u)^{1-p} \mathcal{R}_1^*\left(\left|K_j\cap \cdot\right|^{n-1}\right)(u) d u\nonumber\\
&\geq&\frac{2}{n\omega_n}\int_{S^{n-1}} \left(u\cdot\nu_{K_j}\left(\rho_{K_j}(u)u\right)\right)^{-p}\rho_{K_j}(u)^{1-p}\rho_{K_j}(u)^{n-1} d u\nonumber\\
&=&\frac{2}{n\omega_n}\int_{S^{n-1}} \left(u\cdot\nu_{K_j}\left(\rho_{K_j}(u)u\right)\right)^{-p}\rho_{K_j}(u)^{n-p}d u\nonumber\\
&=&\frac{2}{\omega_n}\widetilde{C}_{p,n}\left(K_j, S^{n-1}\right),\nonumber
\end{eqnarray}
which implies that $\widetilde{C}_{p,n}\left(K_j, S^{n-1}\right)$ stays bounded if $\widetilde{C}^a_{p,1}\left(K_j, S^{n-1}\right)$ stays bounded. From Lemma \ref{BF19} for $\widetilde{C}_{p,n}(\cdot)$, we obtain $K \in \mathcal{K}_{(o)}^n$.
\end{proof}

\section{Solution to the $L_p$ Minkowski problem on ADQ for non-symmetric discrete measures}\label{SSNSDm}

In this section, we solve the existence part of the $L_p$ Minkowski problem on ADQ for non-symmetric discrete measures. The main method follows along the same lines as that of \cite{BF19}. The difference is that we always have to consider $\mathcal{R}_m^*\left(\left|P(z) \cap \cdot\right|^{n-1}\right)$ for polytope $P(z)$ defined in (\ref{4.5}) in proving $o\in\operatorname{int}P(z_0)$, while only $\rho_{P(z)}$ is considered in \cite{BF19}.
%

\begin{theorem}\label{T5.1}
 Let $p>1$ and $m=1,\dots,n-1$, and let $\mu$ be a discrete measure on $S^{n-1}$ that is not concentrated on any closed hemisphere. Then there exists a polytope $P \in \mathcal{K}_{(o)}^n$ such that
 \begin{equation}\label{PsiP1}
 \widetilde{\Psi}_m(P)^{-1} \widetilde{C}^a_{p,m}(P, \cdot)=\mu.
 \end{equation}
 \end{theorem}

We observe that if $P \in \mathcal{K}_o^n$ is a polytope with $ \operatorname{int}P \neq \emptyset$, and $v_1, \ldots, v_l \in S^{n-1}$ are the exterior normals of the facets of $P$ not containing the origin, then
\begin{equation}
\operatorname{supp} \widetilde{C}^a_{p,m}(P, \cdot) =\left\{v_1, \ldots, v_l\right\},\nonumber
\end{equation}
and for $i=1, \ldots, l$,
\begin{equation}
\widetilde{C}^a_{p,m}\left(P,\left\{v_i\right\}\right)= \int_{\langle\cdot\rangle_{P}\left(F\left(P, v_i\right)\right)}\left(u\cdot\nu_P\left(\rho_P(u)u\right)\right)^{-p} \rho_P^{m-p}(u) \mathcal{R}_m^*\left(\left|P\cap \cdot\right|^{n-1}\right)(u) du. \nonumber
\end{equation}

Let $p>1$, $m=1,\dots,n-1$ and $\mu$ be a discrete measure on $S^{n-1}$ that is not concentrated on any closed hemisphere. Let $\operatorname{supp} \mu=\left\{u_1, \ldots, u_k\right\}$, and let $\mu\left(\left\{u_i\right\}\right)=\alpha_i>0$, $i=1, \ldots, k$. For $z=\left(t_1, \ldots, t_k\right) \in (\mathbb{R}_{\geq 0})^k$, we define
\begin{eqnarray}
\mathcal{F}(z)&=&\sum_{i=1}^k \alpha_i t_i^p \label{4.5}\\
P(z)&=&\left\{x \in \mathbb{R}^n: x\cdot u_i \leq t_i,\; \forall i=1, \ldots, k\right\}\label{4.6} \\
Z\;\;\;&=&\left\{z \in\left(\mathbb{R}_{\geq 0}\right)^k: \mathcal{F}(z)=1\right\}.\label{4.7}
\end{eqnarray}

Since $\alpha_i>0$ for $i=1, \ldots, k$, the set $Z$ is compact. Therefore, by the continuity of $\widetilde{\Psi}_m(\cdot)$, see Lemma \ref{L4.1}, there exists $z_0 \in Z$ such that
\begin{equation}\label{4.17}
\widetilde{\Psi}_m\left(P(z_0)\right)=\max \left\{\widetilde{\Psi}_m\left(P(z)\right): z \in Z\right\}.
\end{equation}
We will prove that $o \in \operatorname{int} P\left(z_0\right)$ and there exists a $\lambda_0>0$ such that
\begin{equation}
\widetilde{\Psi}_m\left(\lambda_0 P\left(z_0\right)\right)^{-1} \widetilde{C}^a_{p,m}\left(\lambda_0 P\left(z_0\right), \cdot\right)=\mu.\nonumber
\end{equation}

\begin{lemma}\label{LL6.1}
Let $n>m\geq 2$. If $z_0$ satisfies (\ref{4.17}),  then there exists a constant $c_1$ depending on $n,m,z_0$ such that for  {all} $u\in S^{n-1}$, we have
\begin{equation}\label{Eint}
\mathcal{R}_m^*\left(\left|P(z_0) \cap \cdot\right|^{n-1}\right)(u)\geq c_1.
\end{equation}
\end{lemma}
\begin{proof}
Since $z_0$ satisfies (\ref{4.17}), the interior of $P\left(z_0\right)$ {is nonempty}. Thus, there exists a ball $B^n(x_0,r_0)\subset P(z_0)$, where $r_0$ depends on $z_0$. Moreover, there exists a $R_0>0$ such that $P(z_0)\subset B^n(o,R_0)$. By  Lemma \ref{L5.6}, there exists a constant $c_1$ depending on $n,m,z_0$ such that the inequality (\ref{Eint}) is established for  {all} $u\in S^{n-1}$.
\end{proof}

\begin{lemma}\label{L4.3}
Let $p>1$ and $m=1,\dots,n-1$. If $z_0$ satisfies (\ref{4.17}),  then $o \in \operatorname{int} P\left(z_0\right)$.
\end{lemma}

\begin{proof}
It is obvious from the construction of $Z$ and (\ref{4.6}) that $o \in P\left(z_0\right)$. We will prove this conclusion by contradiction. Assume that $o \in \partial P\left(z_0\right)$. Without loss of generality, we may assume that $z_0=\left(t_1, \ldots, t_k\right) \in\left(\mathbb{R}_{\geq 0}\right)^k$, where there exists $1 \leq i_0<k$ such that $t_1=\cdots=t_{i_0}=0$ and $t_{i_0+1}, \ldots, t_k>0$. For sufficiently small $t>0$, we define
\begin{eqnarray}
&&\alpha=\frac{\alpha_1+\ldots+\alpha_{i_0}}{\alpha_{{i_0}+1}+\ldots+\alpha_k},\notag\\
&& \tilde{z}_t=\left(\overbrace{0, \ldots, 0}^{i_0},\left(t_{{i_0}+1}^p-\alpha t^p\right)^{\frac{1}{p}}, \ldots,\left(t_k^p-\alpha t^p\right)^{\frac{1}{p}}\right),\label{ztilde} \\
&&z_t=\left(\overbrace{t, \ldots, t}^{i_0},\left(t_{{i_0}+1}^p-\alpha t^p\right)^{\frac{1}{p}}, \ldots,\left(t_k^p-\alpha t^p\right)^{\frac{1}{p}}\right).\label{zt}
\end{eqnarray}
Simple substitution shows that $\mathcal{F}(z_t)=1$, so $z_t\in Z$.

We prove that there exist $\tilde{t}_0, \tilde{c}_1, \tilde{c}_2>0$ depending on $n,m,p,\mu$ and $z_0$ such that if $t \in\left(0, \tilde{t}_0\right]$ then
\begin{equation}
\widetilde{\Psi}_m\left(P(\tilde{z}_t)\right)\geq \widetilde{\Psi}_m\left(P(z_0)\right)-\tilde{c}_1 t^p\label{41}
\end{equation}
and
\begin{equation}
\widetilde{\Psi}_m\left(P(z_t)\right)\geq \widetilde{\Psi}_m\left(P(\tilde{z}_t)\right)+\tilde{c}_2 t,\label{42}.
\end{equation}
If (\ref{41}) and (\ref{42}) are established, then
\begin{equation}\label{43}
\widetilde{\Psi}_m \left(P(z_t)\right)\geq \widetilde{\Psi}_m\left(P(z_0)\right)-\tilde{c}_1 t^p+\tilde{c}_2 t,
\end{equation}
 which contradicts the maximality of $\widetilde{\Psi}_m\left(P(z_0)\right)$ when $t\rightarrow 0^+$.

First, we prove (\ref{41}). We choose $R>0$ such that $P\left(z_0\right) \subset \operatorname{int} R B^n$ and $R \geq \max \left\{t_{i_0+1}, \ldots, t_k\right\}$.
 Let $\rho_0=\min \left\{t_{i_0+1}, \ldots, t_k\right\}$.
It follows from $p>1$ and $\rho_0 \leq t_i \leq R$  that there exists $s_0>0$, depending on $z_0, \mu$ and $p$ such that if $t \in\left(0, s_0\right)$, then for $i=i_0+1, \ldots, k$,
\begin{eqnarray}\label{45}
\left(t_i^p-\alpha t^p\right)^{\frac{1}{p}}>t_i\left(1-\alpha\frac{t^p}{t_i^p}\right)=t_i-\alpha t_i^{1-p}t^p>t_i-\alpha\rho_0^{1-p} t^p>\rho_0 / 2.
\end{eqnarray}
 By (\ref{NKo}), it is easily seen that
$$N\left(P\left(z_0\right), o\right)^*=\left\{x \in \mathbb{R}^n: x\cdot u_i\leq 0,\; \forall i=1, \ldots, i_0\right\},$$
and $\rho_{P\left(z_0\right)}(u)>0$ if and only if $u \in N\left(P\left(z_0\right), o\right)^*$.
 It follows from (\ref{ztilde}) and (\ref{45}) that for $t \in\left(0, s_0\right)$,
 \begin{equation}\label{7.13}
\rho_{P\left(\tilde{z}_t\right)}(u)>\rho_0 / 2\;\;{\rm for}\;u\in S^{n-1}\quad {\rm if\;and\;only\;if}\quad u \in N\left(P\left(z_0\right), o\right)^*.
\end{equation}

Let $u \in N\left(P\left(z_0\right), o\right)^* \cap S^{n-1}$. Then $\rho_{P\left(\tilde{z}_t\right)}(u) u$ lies in a facet $F\left(P\left(\tilde{z}_t\right), u_i\right)$ for an $i \in$ $\{i_0+1, \ldots, k\}$ by (\ref{7.13}). Therefore, by (\ref{ztilde}), (\ref{45}), for $t \in\left(0, s_0\right)$,
\begin{equation}\label{Ptildez}
\rho_{P\left(\tilde{z}_t\right)}(u) u\cdot u_i=\left(t_i^p-\alpha t^p\right)^{\frac{1}{p}}>t_i-\alpha\rho_0^{1-p} t^p>\rho_0/2.
\end{equation}

By the above inequality and $\rho_{P\left(\tilde{z}_t\right)}(u) \leq R$, we deduce that $u\cdot u_i\geq \frac{\rho_0}{2 R}$.
Let $s>0$ be defined by $s u\cdot u_i=t_i$. Then $s \geq \rho_{P\left(z_0\right)}(u)$, and hence by (\ref{Ptildez}),
$$
s-\rho_{P\left(\tilde{z}_t\right)}(u)=\frac{s u\cdot u_i-\rho_{P\left(\tilde{z}_t\right)}(u) u\cdot u_i}{u\cdot u_i} \leq \frac{t_i-\left(t_i-\alpha\rho_0^{1-p} t^p\right)}{u\cdot u_i} \leq 2 R \alpha\rho^{-p}_0 t^p.
$$
Thus,
\begin{equation}\label{Ptz}
\rho_{P\left(\tilde{z}_t\right)}(u) \geq \rho_{P\left(z_0\right)}(u)-2 R \alpha\rho^{-p}_0 t^p.
\end{equation}

We choose $t_0>0$ with $t_0 \leq s_0$ depending on $z_0$ and $p$ such that
\begin{equation}\label{E65}
2m R \alpha\rho^{-p}_0t_0^p<\rho_0 / 2.
\end{equation}

If $t \in\left(0, t_0\right)$ and $u \in N\left(P\left(z_0\right), o\right)^* \cap S^{n-1}$, then by (\ref{Ptz}), (\ref{E65}) and Bernoulli's inequality,
\begin{eqnarray}\label{Evar}
\rho_{P\left(\tilde{z}_t\right)}(u)^m \geq\left(\rho_{P\left(z_0\right)}(u)-2 R \alpha\rho^{-p}_0t^p\right)^m \geq\rho_{P\left(z_0\right)}(u)^m-2m R^m \alpha\rho^{-p}_0 t^p.
\end{eqnarray}

If $u\notin N\left(P\left(z_0\right), o\right)^*$, then $\rho_{P\left(\tilde{z}_t\right)}(u)=0=\rho_{P\left(z_0\right)}(u)$.
Thus, for $\xi\in G(n,m)$, by (\ref{Evar}),
\begin{eqnarray}\label{EPleft}
\left|P\left(\tilde{z}_t\right)\cap \xi\right|
\geq\left|P\left(z_0\right)\cap \xi\right|-2 R^m \alpha\rho^{-p}_0\mathcal{H}^{m-1}(\xi\cap S^{n-1}\cap N\left(P\left(z_0\right), o\right)^*)t^p.
\end{eqnarray}

Since $\rho_0 \leq \rho_{P\left(z_0\right)}(u) \leq R$ for $u \in N\left(P\left(z_0\right), o\right)^* \cap S^{n-1}$, then for $\xi\in G(n,m)$,
\begin{eqnarray}\label{E82}
\left|P\left(z_0\right)\cap \xi\right|\geq\frac{1}{m}\rho_0^m\mathcal{H}^{m-1}\left(\xi\cap S^{n-1}\cap N\left(P\left(z_0\right), o\right)^*\right).
\end{eqnarray}

By (\ref{EPleft}) and (\ref{E82}), there exists a sufficiently small $t_1<t_0$ such that for $t\in (0,t_1)$ and any $\xi\in G(n,m)$ and $|\xi\cap P(z_0)|>0$, we have
\begin{equation}\label{E69}
 \left|P\left(z_0\right)\cap \xi\right|-2 R^m \alpha\rho^{-p}_0\mathcal{H}^{m-1}(\xi\cap S^{n-1}\cap N\left(P\left(z_0\right), o\right)^*)t^p>0.
\end{equation}

By (\ref{EPleft}) and (\ref{E69}), if  $\xi\in G(n,m)$ and $\left|P\left(z_0\right)\cap \xi\right|>0$, then
\begin{eqnarray}\label{E70}
\left|P\left(\tilde{z}_t\right)\cap \xi\right|^n &\geq& \left(\left|P\left(z_0\right)\cap \xi\right|-2R^m \alpha\rho^{-p}_0\mathcal{H}^{m-1}(\xi\cap S^{n-1}\cap N\left(P\left(z_0\right), o\right)^*)t^p\right)^n \nonumber\\
&\geq& \left|P\left(z_0\right)\cap \xi\right|^n-2n R^m \alpha\rho^{-p}_0\left|P\left(z_0\right)\cap \xi\right|^{n-1}\mathcal{H}^{m-1}(\xi\cap S^{n-1}\cap N\left(P\left(z_0\right), o\right)^*)t^p\nonumber\\
&\geq& \left|P\left(z_0\right)\cap \xi\right|^n-2m^{1-n}n R^{mn} \alpha\rho^{-p}_0\mathcal{H}^{m-1}(\xi\cap S^{n-1}\cap N\left(P\left(z_0\right), o\right)^*)^n t^p.
\end{eqnarray}
Let
\begin{equation}\label{E71}
\tilde{c}_1=2m^{1-n}n R^{mn} \alpha\rho^{-p}_0\int_{G(n,m)}\mathcal{H}^{m-1}(\xi\cap S^{n-1}\cap N\left(P\left(z_0\right). o\right)^*)^nd\nu_m(\xi).
\end{equation}
By (\ref{E70}) and (\ref{E71}), we have
\begin{equation}
\widetilde{\Psi}_m\left(\tilde{z}_t\right)=\int_{G(n,m)}\left|P\left(\tilde{z}_t\right)\cap \xi\right|^nd\nu_m(\xi)\geq \int_{G(n,m)}\left|P\left(z_0\right)\cap \xi\right|^nd\nu_m(\xi)-\tilde{c}_1\cdot t^p=\widetilde{\Psi}_m\left(z_0\right)-\tilde{c}_1\cdot t^p,\nonumber
\end{equation}
which is the desired inequality (\ref{41}).

The proof of (\ref{42}) is based on the idea of constructing a set $\widetilde{U}_t \subset S^{n-1}$ for sufficiently small $t>0$ whose $\mathcal{H}^{n-1}$ measure is of order $t$, and if $u \in \widetilde{U}_t$, then $\rho_{P\left(z_t\right)}(u) \geq r$ for a suitable constant $r>0$ depending $n,m,p,\mu$ and $P\left(z_0\right)$ while $\rho_{P\left(\tilde{z}_t\right)}(u)=0$.
Without loss of generality we can assume that $\operatorname{dim} F\left(P\left(z_0\right), u_1\right)=n-1$. In particular, there exist $r>0$ and $y_0 \in F\left(P\left(z_0\right), u_1\right) \backslash\{o\}$ such that
\begin{equation}\label{8r}
y_0\cdot u_i\leq h_{P\left(z_0\right)}\left(u_i\right)-8 r \text { for } i=2, \ldots, k.
\end{equation}
For $v=y_0 /\left\|y_0\right\| \in S^{n-1} \cap u_1^{\perp}$, we consider $y=y_0+4 r v$, and hence $4 r \leq\|y\| \leq R$, and by (\ref{8r}),
\begin{equation}\label{4r}
y\cdot u_i\leq h_{P\left(z_0\right)}\left(u_i\right)-4 r \quad \text { for } i=2, \ldots, k.
\end{equation}
Note that $P\left(\tilde{z}_t\right) \rightarrow P\left(z_0\right)$ as $t \rightarrow 0^{+}$and also $P\left(\tilde{z}_t\right) \subset P\left(z_0\right)$ for $t>0$. Therefore there exists a positive $t_1 \leq \min \left\{r, t_0\right\}$, depending only on $p, m, \mu$ and $z_0$ such that if $t \in\left(0, t_1\right]$, then by (\ref{4r}),
\begin{equation}\label{46}
y\cdot u_i\leq h_{P\left(\tilde{z}_t\right)}\left(u_i\right)-2 r \quad \text { for } i=2, \ldots, k \quad \text { and } P\left(z_t\right) \subset R B^n.
\end{equation}

 Let the $(n-2)$-dimensional unit ball $U$ be defined as $U=u_1^{\perp} \cap v^{\perp} \cap B^n$.
Then we have that $y+r U \subset F\left(P\left(z_0\right), u_1\right)$ and $(y+r U)+r\left[o, u_1\right] \subset y+2 r B^n$. Let
 $U_t=(y+r U)+t\left(o, u_1\right]$
 be the $(n-1)$-dimensional right spherical cylinder of height $t<\min \left\{t_1, r\right\}$, whose base $y+r U$ does not belong to $U_t$. We deduce from (\ref{46}) and $h_{P\left(z_t\right)}\left(u_1\right)=t$ that
 \begin{equation}\label{Ut}
 U_t \subset P\left(z_t\right) \backslash N\left(P\left(z_0\right), o\right)^* \subset P\left(z_t\right) \backslash P\left(\tilde{z}_t\right).
 \end{equation}

Let $\widetilde{U}_t$ be the radial projection of $U_t$ to $S^{n-1}$. For $x \in U_t$, we have $x\cdot v= {|y|} \geq 4 r$ and $ {|x|} \leq R$, therefore
\begin{equation}\label{wUt}
\mathcal{H}^{n-1}\left(\widetilde{U}_t\right)=\int_{U_t} x\cdot v {|x|}^{-n} d \mathcal{H}^{n-1}(x) \geq \frac{4 r \mathcal{H}^{n-1}\left(U_t\right)}{R^n}=\frac{4 r^{n-1} \omega_{n-2}}{R^n} t.
\end{equation}
By (\ref{Ut}), if $u \in \widetilde{U}_t$, then
\begin{equation}\label{Pzt}
\rho_{P\left(\tilde{z}_t\right)}(u)=0\;{\rm and}\;\rho_{P\left(z_t\right)}(u) \geq\|y\| \geq 4 r.
\end{equation}

When $n>m\geq 2$, by (\ref{E4.4}), (\ref{Pzt}), (\ref{wUt}) and Lemma \ref{LL6.1},  there exists a constant $c_1>0$ such that
\begin{eqnarray}\label{E7.23}
\widetilde{\Psi}_m\left(P\left(z_t\right)\right)
&=&\frac{1}{m} \int_{S^{n-1} \backslash \widetilde{U}_t} \rho_{P\left(z_t\right)}^m(u)\mathcal{R}_m^*\left(\left|P(z_t) \cap \cdot\right|^{n-1}\right)(u)   d \mathcal{H}^{n-1}(u)\nonumber\\
&&+\frac{1}{m} \int_{\widetilde{U}_t} \rho_{P\left(z_t\right)}^m(u) \mathcal{R}_m^*\left(\left|P(z_t) \cap \cdot\right|^{n-1}\right)(u) d \mathcal{H}^{n-1}(u)\nonumber\\
&\geq&\frac{1}{m} \int_{S^{n-1}} \rho_{P\left(\tilde{z}_t\right)}^m(u) \mathcal{R}_m^*\left(\left|P(\tilde{z}_t) \cap \cdot\right|^{n-1}\right)(u)d \mathcal{H}^{n-1}(u)\nonumber\\
&&+\frac{1}{m} \int_{\widetilde{U}_t} \rho_{P\left(z_t\right)}^m(u)\mathcal{R}_m^*\left(\left|P(\tilde{z}_t) \cap \cdot\right|^{n-1}\right)(u)d \mathcal{H}^{n-1}(u)\nonumber\\
&\geq&\widetilde{\Psi}_m\left(P\left(\tilde{z}_t\right)\right)+\frac{4^{m+1} r^{m+n-1} \omega_{n-2}c_1}{mR^n}t.
\end{eqnarray}

 When $n\geq 2$ and $m=1$, by (\ref{Pzt}) and (\ref{wUt}),
 \begin{eqnarray}\label{m1}
 \frac{1}{m} \int_{\widetilde{U}_t} \rho_{P\left(z_t\right)}^m(u) \mathcal{R}_1^*\left(\left|P(z_t) \cap \cdot\right|^{n-1}\right)(u) d \mathcal{H}^{n-1}(u)\nonumber\\
\geq\frac{2}{n\omega_n}\int_{\widetilde{U}_t} \rho_{P\left(z_t\right)}^{n}(u) d \mathcal{H}^{n-1}(u)\geq\frac{2^{2n+3}r^{2n-1} \omega_{n-2}}{n\omega_nR^n} t.
 \end{eqnarray}
Let
$$\tilde{c}_2=\min\left\{\frac{4^{m+1} r^{m+n-1} \omega_{n-2}c_1}{mR^n},\frac{2^{2n+3}r^{2n-1} \omega_{n-2}}{n\omega_nR^n}\right\}.$$
By (\ref{E7.23}) and (\ref{m1}), we get the desired inequality (\ref{42}).
\end{proof}

As we already know that $o \in \operatorname{int} P\left(z_0\right)$ by Lemma \ref{L4.3}, we can freely decrease $h_{P\left(z_0\right)}\left(u_i\right)$ for $i=1, \ldots, k$, and increase it if $\operatorname{dim} F\left(P\left(z_0\right), u_i\right)=n-1$. To control what happens to $\widetilde{\Psi}_m\left(P(z)\right)$ when we perturb $P\left(z_0\right)$, we use Lemma \ref{L3.3}, which is a consequence of Corollary \ref{C2.3} when $p=1$.
%

\begin{lemma}\label{L3.3} If $m=1,\dots,n-1$, $\delta \in(0,1)$ and $z_t=\left(z_1(t), \ldots, z_k(t)\right) \in$ $\mathbb{R}_{+}^k$ for $t \in(-\delta, \delta)$ are such that $\lim _{t \rightarrow 0} \frac{z_i(t)-z_i(0)}{t}=z_i^{\prime}(0) \in \mathbb{R}$ for $i=1, \ldots, k$ exists, then the $P\left(z_t\right)$ defined in (\ref{4.6}) satisfies that
$$
\lim _{t \rightarrow 0} \frac{\widetilde{\Psi}_m\left(P\left(z_t\right)\right)-\widetilde{\Psi}_m\left(P\left(z_0\right)\right)}{t}=n \sum_{i=1}^k \frac{z_i^{\prime}(0)}{h_{P\left(z_0\right)}\left(u_i\right)} \cdot \widetilde{C}^a_m\left(P\left(z_0\right),\left\{u_i\right\}\right).
$$
\end{lemma}


The following lemma give a stronger conclusion, i.e., the support sets $F\left(P\left(z_0\right),\cdot\right)$ with respect to every directions $u_i$, $i=1,\dots,k$, are facets of $P\left(z_0\right)$. This shows that $\widetilde{C}^a_m\left(P\left(z_0\right),\left\{u_i\right\}\right)>0$ for every $u_i$, $i=1,\dots,k$.

\begin{lemma} \label{L4.5}
If $p>1$ and $m=1,\dots,n-1$, then $\operatorname{dim} F\left(P\left(z_0\right), u_i\right)=n-1$  for $i=1, \ldots, k$.
\end{lemma}
\begin{proof}
We assume that $\operatorname{dim} F\left(P\left(z_0\right), u_1\right)<n-1$, and obtain a contradiction. Let $z_0=\left(t_1, \ldots, t_k\right) \in \mathbb{R}_+^k$.  We may assume that $\operatorname{dim} F\left(P\left(z_0\right), u_k\right)=n-1$. For small $t \geq 0$, we consider $\tilde{z}(t)=\left(t_1-t, t_2, \ldots, t_k\right)$,
and $\theta(t)=\mathcal{F}\left(\tilde{z}(t)\right)$.
In particular, $\theta(0)=1$ and
$\theta^{\prime}(0)=-p \alpha_1 t_1^{p-1}$. Let
$z(t)=\theta(t)^{-1 / p} \tilde{z}(t)=\left(z_1(t), \ldots, z_k(t)\right)$.
Then $z(t)\in Z$, $\left.\frac{d}{d t} \theta(t)^{-1 / p}\right|_{t=0^{+}}=\alpha_1 t_1^{p-1}$ and $z_i^{\prime}(0)=\alpha_1 t_1^{p-1} t_i>0$ for $i=2, \ldots, k$. By Lemma \ref{L3.3} and $\widetilde{C}^a_m\left(P\left(z_0\right),\left\{u_1\right\}\right)=0$,  we conclude that
\begin{eqnarray}
&&\lim _{t \rightarrow 0^{+}} \frac{\widetilde{\Psi}_m(P(z(t)))-\widetilde{\Psi}_m\left(P\left(z_0\right)\right)}{t}
=n\sum_{i=2}^k \frac{z_i^{\prime}(0)}{h_{P\left(z_0\right)}\left(u_i\right)} \cdot \widetilde{C}^a_m\left(P\left(z_0\right),\left\{u_i\right\}\right)\nonumber\\
&&~~~~~~~~~~~~\geq\frac{nz_k^{\prime}(0)}{h_{P\left(z_0\right)}\left(u_k\right)} \cdot \widetilde{C}^a_m\left(P\left(z_0\right),\left\{u_k\right\}\right)
=n\alpha_1 t_1^{p-1}\cdot \widetilde{C}^a_m\left(P\left(z_0\right),\left\{u_k\right\}\right)>0,\nonumber
\end{eqnarray}
which implies $\widetilde{\Psi}_m(P(z(t)))>\widetilde{\Psi}_m\left(P\left(z_0\right)\right)$ when $t>0$ sufficiently small.  This leads to a contradiction with the optimality of $z_0$.
\end{proof}

\noindent {\bf Proof of Theorem \ref{T5.1}.} On account of Lemmas \ref{L4.3} and \ref{L4.5}, we have $\operatorname{dim} F\left(P\left(z_0\right), u_i\right)=n-1$ for $i=1, \ldots, k$, $o \in \operatorname{int} P\left(z_0\right)$. Let $\left(s_1, \ldots, s_k\right) \in \mathbb{R}^k$ satisfy $\sum_{i=1}^k s_i \alpha_i t_i^{p-1}=0$ and  not all $s_i$ are zero. If $t \in(-\varepsilon, \varepsilon)$ for small $\varepsilon>0$, then consider $\tilde{z}(t)=\left(t_1+s_1 t, \ldots, t_k+s_k t\right)$
and $\theta(t)=\mathcal{F}(\tilde{z}(t))$. In particular, $\theta(0)=1$ and
$\theta^{\prime}(0)=p \sum_{i=1}^k s_i \alpha_i t_i^{p-1}=0$.
Let $z(t)=\theta(t)^{-1 / p} \tilde{z}(t)=\left(z_1(t), \ldots, z_k(t)\right) \in Z$.
Then $z_i^{\prime}(0)=s_i$ for $i=1, \ldots, k$. By Lemma \ref{L3.3} and the optimality of $z_0$, we that
\begin{equation}\label{47}
\lim _{t \rightarrow 0} \frac{\tilde{\Psi}_m(P(z(t)))-\widetilde{\Psi}_m\left(P\left(z_0\right)\right)}{t}=n\sum_{i=1}^k \frac{s_i}{t_i} \widetilde{C}^a_m\left(P\left(z_0\right),\left\{u_i\right\}\right)=0.
\end{equation}

In particular, (\ref{47}) holds whenever $\left(s_1, \ldots, s_k\right) \in \mathbb{R}^k \backslash\{o\}$ satisfies $\sum_{i=1}^k s_i \alpha_i t_i^{p-1}=0$. Therefore, there exists a $\lambda \in \mathbb{R}$ such that
\begin{equation}\label{lambdacdot}
\lambda \widetilde{C}^a_m\left(P\left(z_0\right),\left\{u_i\right\}\right)=\alpha_i t_i^{p} \quad for\quad i=1, \ldots, k.
\end{equation}

 Since $\lambda>0$ and $p>1$, there exists a $\lambda_0>0$ such that $\lambda=\lambda_0^{-p} \widetilde{\Psi}_m\left(P\left(z_0\right)\right)^{-1}$. Hence by (\ref{lambdacdot}), the $mn$-homogeneity of  $\widetilde{C}_{m}^{\mathrm{a}}$ and the $mn$-homogeneity of  $\widetilde{\Psi}_m$, we can conclude that
 $$
\alpha_i=\widetilde{\Psi}_m\left(\lambda_0 P\left(z_0\right)\right)^{-1} h_{\lambda_0 P\left(z_0\right)}\left(u_i\right)^{-p} \widetilde{C}^a_m\left(\lambda_0 P\left(z_0\right),\left\{u_i\right\}\right) \quad \text { for } i=1, \ldots, k.
$$
This implies that (\ref{PsiP1}) is established if let $P=\lambda_0 P(z_0)$.
\qed
$\;$\\

\noindent{\bf Proof of Theorem \ref{T1.1}.} Let $p \neq mn$. According to Theorem \ref{T5.1}, there exists a polytope $P_0 \in \mathcal{K}_{(o)}^n$ such that $\widetilde{\Psi}_m\left(P_0\right)^{-1} \widetilde{C}^a_{p,m}\left(P_0, \cdot\right)=\mu$. Let $\lambda^{mn-p}=\widetilde{\Psi}_m\left(P_0\right)^{-1}$ and $P=\lambda P_0$. Then by the $(mn-p)$-homogeneity of  $\widetilde{C}_{p,m}^{\mathrm{a}}$, we have
$\widetilde{C}^a_{p,m}(P, \cdot)=\lambda^{mn-p} \widetilde{C}^a_{p,m}\left(P_0, \cdot\right)=\mu$.
\qed\\

In \cite[Theorem 8.3]{LYZ18},  Lutwak-Yang-Zhang established the uniqueness of the solution on the discrete $L_p$ dual Minkowski problem when $p>q$. We  establish  the uniqueness of the solution to the $L_p$ Minkowski problem on ADQ for the case of polytopes when $p>mn$.
\begin{theorem}
Let $P, P^{\prime} \in \mathcal{K}_o^n$ be polytopes containing the origin in their interiors. Suppose
\begin{equation*}
\widetilde{C}^a_{p,m}(P,\cdot)=\widetilde{C}^a_{p,m}\left(P^{\prime},\cdot\right).
\end{equation*}
Then $P=P^{\prime}$ when $p>mn$ and $P^{\prime}$ is a dilate of $P$ when $p=mn$.
\end{theorem}
\begin{proof}
By (\ref{E3.12}), $\widetilde{C}^a_{p,m}(P,\cdot)=\widetilde{C}^a_{p,m}\left(P^{\prime},\cdot\right)$ means that $P$ and $P^{\prime}$ must have the same outer unit normals. Let $u_1, \ldots, u_k$ be the outer unit normals of $P$ and $P^{\prime}$. Then
$$
\widetilde{C}^a_{p,m}(P,\cdot)=\widetilde{C}^a_{p,m}\left(P^{\prime},\cdot\right)=\sum_{i=1}^k \alpha_i \delta_{u_i},
$$
where
\begin{eqnarray}\label{E7.35}
\alpha_i&=&h_P\left(u_i\right)^{-p} \int_{S^{n-1} \cap \Delta_i} \rho_P^m(u)\mathcal{R}_m^*\left(\left|P \cap \cdot\right|^{n-1}\right)(u) d u\nonumber\\
&=&h_{P^{\prime}}\left(u_i\right)^{-p} \int_{S^{n-1} \cap \Delta_i^{\prime}} \rho_{P^{\prime}}^m(u)  \mathcal{R}_m^*\left(\left|P \cap \cdot\right|^{n-1}\right)(u) d u
\end{eqnarray}
and  $\delta_{u_i}$ denotes the delta measure concentrated at $u_i$, and $\Delta_i$ and $\Delta_i^{\prime}$ are the cones formed by the origin and the facets of $P$ and $P^{\prime}$ with normal $u_i$, respectively.

Assume that $P \neq P^{\prime}$. Then $P\subset P^{\prime}$ is not possible. Let $\lambda$ be the maximal number so that $\lambda P \subseteq P^{\prime}$. Then $\lambda<1$. Since $\lambda P$ and $P^{\prime}$ have the same outer unit normals, there is a facet of $\lambda P$ that is contained in a facet of $P^{\prime}$. Denote the outer unit normal of those facets by $u_{i_1}$. We have
$h_{\lambda P}\left(u_{i_1}\right)=h_{P^{\prime}}\left(u_{i_1}\right)$, $\Delta_{i_1}\subseteq \Delta_{i_1}^{\prime}$
 and
$$\rho_{\lambda P}(u)=\rho_{P^{\prime}}(u)\quad {\rm for\;all}\quad u \in \Delta_{i_1}.$$
Moreover, by $\lambda P\subseteq P^{\prime}$, we have
$$\mathcal{R}_m^*\left(\left|\lambda P \cap \cdot\right|^{n-1}\right)(u)\leq\mathcal{R}_m^*\left(\left|P^{\prime} \cap \cdot\right|^{n-1}\right)(u).$$
Therefore,
\begin{eqnarray}\label{E7.36}
&&h_{\lambda P}\left(u_{i_1}\right)^{-p} \int_{S^{n-1} \cap \Delta_{i_1}} \rho_{\lambda P}^m(u) \mathcal{R}_m^*\left(\left|\lambda P \cap \cdot\right|^{n-1}\right)(u) d u\nonumber\\
& \leq& h_{P^{\prime}}\left(u_{i_1}\right)^{-p} \int_{S^{n-1} \cap \Delta_{i_1}^{\prime}} \rho_{P^{\prime}}^m(u) \mathcal{R}_m^*\left(\left|P^{\prime} \cap \cdot\right|^{n-1}\right)(u) d u,
\end{eqnarray}
with equality only if $\Delta_{i_1}=\Delta_{i_1}^{\prime}$. Combining (\ref{E7.36}) with (\ref{E7.35})  yields $\lambda^{mn-p}\leq 1$. If $p>mn$, then by $\lambda<1$,
we have that $\lambda^{mn-p}>1$. This is a contradiction.

If $p=mn$, the (\ref{E7.35}) gives equality in (\ref{E7.36}). Thus,
$\Delta_{i_1}=\Delta^{\prime}_{i_1}$, which implies that the facets of
$\lambda P$ and $P^{\prime}$ with outer unit normal $u_{i_1}$
are the same. Let $u_{i_2}$ be the outer unit normal to
a facet adjacent to the facet whose outer unit normal is
$u_{i_1}$. Since $u_{i_2}$ is the same outer unit normal of $\lambda P$ and $P^{\prime}$ and $F\left(\lambda P,u_{i_1}\right)= F\left(P^{\prime},u_{i_1}\right)$, we have
$F\left(\lambda P,u_{i_2}\right)\subset F\left(P^{\prime},u_{i_2}\right)$. The same
argument yields that these two facets are also the same. Repeating the previous argument, we get $\lambda P=P^{\prime}$.
\end{proof}

\section{Solution to the $L_p$ Minkowski problem on ADQ for non-symmetric general measures}\label{SSNSGm}

For $w \in S^{n-1}$ and $t \in(-1,1)$, we write

$$
\Omega(w, t)=\left\{u \in S^{n-1}:\; u\cdot w>t\right\}.
$$

\begin{lemma}\label{L6.1}
 For a finite Borel measure $\mu$ on $S^{n-1}$ not concentrated on a closed hemi-sphere, there exists $t \in(0,1)$ such that  $\mu(\Omega(w, t))>t$ for {all} $w \in S^{n-1}$.
\end{lemma}
\begin{proof}
Assume that there exist a sequence of $\{w_i\}_{i=1}^{\infty}\subset S^{n-1}$ such that
\begin{equation}\label{lim1}
\mu(\Omega(w_i, 1/i))\leq 1/i.
\end{equation}
Then by $w_i\in S^{n-1}$, there exists a subsequence $\{w_{i_j}\}_{j=1}^{\infty}$ of $\{w_i\}_{i=1}^{\infty}\subset S^{n-1}$ and $w_0\in S^{n-1}$ such that
$\lim_{j\rightarrow\infty}w_{i_j}=w_0$.
By (\ref{lim1}), we have $\mu(\Omega(w_0, 0))\leq 0$ when $j\rightarrow\infty$.
This contradicts  the fact that $\mu$ is not concentrated on a closed hemisphere.
\end{proof}

%

\begin{lemma}\label{L7.2}
 For $p>n+m(n-2)$, $m=1,2,\dots,n-1$, and finite Borel measure $\mu$ on $S^{n-1}$ not concentrated on a closed hemisphere, there exists a convex body $K \in \mathcal{K}_{(o)}^n$ such that
\begin{equation}\label{E7.2}
d \mu=\widetilde{\Psi}_m(K)^{-1} d \widetilde{C}^a_{p,m}(K, \cdot).
\end{equation}
\end{lemma}

\begin{proof}
 Let $\mu_j$, $j\in\mathbb{N}$, be a sequence of discrete measures converging to $\mu$ that are not concentrated on any closed hemispheres. It follows from Theorem \ref{T5.1} that there exists a polytope $P_j \in \mathcal{K}_{(o)}^n$ such that
\begin{equation}\label{E6.1}
d \mu_j=\widetilde{\Psi}_m\left(P_j\right)^{-1} d \widetilde{C}^a_{p, m}\left(P_j, \cdot\right)
\end{equation}
for each $j$, and hence we may assume that
\begin{equation}\label{E6.2}
\widetilde{\Psi}_m\left(P_j\right)^{-1} \widetilde{C}^a_{p,m}\left(P_j, S^{n-1}\right)<2 \mu\left(S^{n-1}\right).
\end{equation}

Next, we prove that $P_j$ is bounded.
Let $R_j=\max _{x \in P_j} {|x|}$. Let $v_j \in S^{n-1}$ such that $R_j v_j \in P_j$. Since $\{v_j\}_{j=1}^{\infty}\subset S^{n-1}$, there exists a subsequence (still denoted by $\{v_j\}$) tends to $v \in S^{n-1}$. By Lemma \ref{L6.1}, there exist $s, t>0$ such that $\mu(\Omega(v, 2 t))>2 s$. As $v_j$ tends to $v \in S^{n-1}$ and $\mu_j$ tends weakly to $\mu$, we have
$\Omega(v, 2 t) \subset \Omega\left(v_j, t\right)\quad {\rm and}\quad\mu_j(\Omega(v, 2 t))>s$.
 Therefore for sufficiently large $j$, we have  $\mu_j\left(\Omega\left(v_j, t\right)\right)>s$.

From $h_{P_j}(u) \geq(u\cdot R_j v_j) \geq R_j t$ for $u \in \Omega\left(v_j, t\right)$, (\ref{E6.1}), (\ref{2.4}) and (\ref{E4.4}), we conclude that
\begin{eqnarray}
s<\mu_j\left(\Omega\left(v_j, t\right)\right)&=&\int_{\Omega\left(v_j, t\right)} \widetilde{\Psi}_m\left(P_j\right)^{-1}h_{P_j}^{-p}(u) d \widetilde{C}^a_m\left(P_j, u\right)\nonumber\\
 &\leq& R_j^{-p} t^{-p} \widetilde{C}^a_m\left(P_j, S^{n-1}\right)\widetilde{\Psi}_m\left(P_j\right)^{-1}=mR_j^{-p} t^{-p}.
\end{eqnarray}
Thus, $R_j^p \leq ms^{-1} t^{-p}$. Thus, $P_j$ is bounded, i.e., there exists $R>0$ such that
$P_j\subset RB^n$.

Since $P_j$ is  bounded, there exists a subsequence of (still denoted by $\{P_j\}$) tends to a compact convex set $K \in \mathcal{K}_o^n$ with $K \subset R B^n$.


We deduce from (\ref{E6.2}) and Lemma \ref{L6.2} that
\begin{equation}\label{8.9}
r(P_j)^{p-n-m(n-2)}\geq c_1 \widetilde{C}^a_{p, m}\left(P_j, S^{n-1}\right)^{-1} \widetilde{\Psi}_m(P_j) \geq \frac{c_1}{2} \mu\left(S^{n-1}\right)^{-1}.
\end{equation}
Since $p>n+m(n-2)$, there exists a constant $c_2>0$ independent of $j$ such that
$|P_j|\geq c_2$.

Since $P_j$ convergence to $K$ with respect to Hausdorff metric and $|P_j|\geq c_2>0$, by the continuity of volume with respect to Hausdorff metric, we have $|K|\geq c_2>0$,
which implies that $\operatorname{int}K\neq \emptyset$.

We observe that $h_{P_j}^p$ tends uniformly to $h_K^p$, and hence also $\widetilde{\Psi}_m\left(P_j\right) h_{P_j}^p$ tends uniformly to $\widetilde{\Psi}_m(K) h_K^p$ by Lemma \ref{L4.1}. Therefore given any continuous function $f$, we have
\begin{equation}\label{E6.3}
\lim _{j \rightarrow \infty} \int_{S^{n-1}} f(u) \widetilde{\Psi}_m\left(P_j\right) h_{P_j}^p(u) d \mu_j=\int_{S^{n-1}} f(u) \widetilde{\Psi}_m(K) h_K^p(u) d \mu.
\end{equation}
Moreover, by (\ref{E6.1}) and Proposition \ref{P2.6}, i.e., the weak convergence of $\widetilde{C}^a_m$, we have
\begin{equation}\label{E6.4}
\lim _{j \rightarrow \infty} \int_{S^{n-1}} f(u) \widetilde{\Psi}_m\left(P_j\right) h_{P_j}^p(u) d \mu_j=\lim _{j \rightarrow \infty} \int_{S^{n-1}} f(u) d \widetilde{C}^a_m(P_j, u)=\int_{S^{n-1}} f(u) d \widetilde{C}^a_m(K, u).
\end{equation}
It follows from (\ref{E6.3}) and (\ref{E6.4}) that
$$
\int_{S^{n-1}} f(u) \widetilde{\Psi}_m(K) h_K^p(u) d \mu=\int_{S^{n-1}} f(u) d \widetilde{C}^a_m(K, u).
$$
Since the last property holds for all continuous function $f$, we conclude that $
\widetilde{\Psi}_m(K) h_K^p d \mu=d \widetilde{C}^a_m(K, \cdot)$. Thus, (\ref{E7.2}) is established by (\ref{E25}).
 Moreover, by Lemma \ref{L6.4}, we have $K \in \mathcal{K}_{(o)}^n$ if $p\geq n+m(n-2)$.
\end{proof}

\noindent{\bf Proof of Theorem \ref{T1.3}} Let $p>n+m(n-2)$, $m=1,\dots,n-1$. According to Lemma \ref{L7.2}, there exists a $K_0 \in \mathcal{K}_{(o)}^n$ such that $\widetilde{\Psi}_m\left(K_0\right)^{-1} \widetilde{C}^a_{p, m}\left(K_0, \cdot\right)=\mu$. Let $$\lambda=\widetilde{\Psi}_m\left(K_0\right)^{\frac{-1}{mn-p}}\;\;{\rm and}\;\;K=\lambda K_0,$$ we have
$$
\widetilde{C}^a_{p,m}(K, \cdot)=\lambda^{mn-p} \widetilde{C}^a_{p,m}\left(K_0, \cdot\right)=\widetilde{\Psi}_m\left(K_0\right)^{-1} \widetilde{C}^a_{p,m}\left(K_0, \cdot\right)=\mu.
$$
\qed

\section{Solution to the $L_p$ Minkowski problem on ADQ for symmetric measures}\label{SSSm}
First, we consider the solution to a maximization problem.
Let $p\geq 0$ and $\mu$ be a non-zero, finite, even Borel measure on $S^{n-1}$ and let $m=1,\dots, n-1$. Define the functional $
J_{m,p,\mu}: \mathcal{K}_{(e)}^n \longrightarrow \mathbb{R}$
by letting
\begin{equation}\label{EJmpmu}
J_{m,p,\mu}(K)=\frac{1}{m n} \log \widetilde{\Psi}_m(K)+E_{p,\mu}(K)
\end{equation}
where
\begin{equation}\label{JEpmu}
E_{p,\mu}(K)=-\frac{1}{p} \log\left(\frac{1}{|\mu|} \int_{S^{n-1}} h^p_K(v) d \mu(v)\right),\;{\rm when}\;p>0
\end{equation}
and
\begin{equation}
E_{p,\mu}(K)=-\frac{1}{|\mu|} \int_{S^{n-1}}\log h_K(v) d \mu(v),\;{\rm when}\;p=0.
\end{equation}

It is obvious that $J_{m,p,\mu}$ is homogenous of degree $0$; i.e., $J_{m,p,\mu}(c K)=J_{m,p,\mu}(K)$ for all real $c>0$. When $K=B^n$, by (\ref{EJmpmu}), (\ref{JEpmu}), (\ref{E4.4}) and (\ref{Rmstar}), we have
\begin{equation}\label{E110}
J_{m,p,\mu}(B^n)=\frac{1}{m}\log\omega_m.
\end{equation}
\noindent {\bf Maximization Problem.} Given a finite even Borel measure $\mu$ on $S^{n-1}$, does there exist a $K_0 \in \mathcal{K}_{(e)}^n$ such that
\begin{equation}\label{4.3}
J_{m,p,\mu}\left(K_0\right)=\sup \left\{J_{m,p,\mu}(K): K \in \mathcal{K}_{(e)}^n\right\} ?
\end{equation}

We are now in a position to show that the solution of the Maximization Problem is the solution of the $L_p$ Minkowski problem on ADQ  for symmetric measures.

\begin{lemma}\label{T4.1} Let $\mu$ be a non-zero, finite, even Borel measure and let $p\geq0$ and $m=1,\dots, n-1$. If there exists a $K_0 \in \mathcal{K}_{(e)}^n$ that is a maximizer of problem (\ref{4.3}),
then there is a convex body $\bar{K} \in \mathcal{K}_{(e)}^n$ such that $\mu=\widetilde{C}_{p,m}^{\mathrm{a}}(\bar{K}, \cdot)$.
\end{lemma}

\begin{proof}
The case of $p=0$ have been proved in the paper \cite{CLWX}, here we only consider the case of $p>0$.
 Let $f \in C_e\left(S^{n-1}\right)$. Let $h_0=h_{K_0}$, and define $h_{t}$ as in (\ref{3.1})
for a suitable chosen $\delta>0$. Let $K_{t}$ be the $L_p$ family of Wulff  shapes defined by (\ref{3.2}).
Consider the function $\phi:(-\delta, \delta) \rightarrow \mathbb{R}$, defined by
$$
\phi(t)=\frac{1}{m n} \log \widetilde{\Psi}_m\left(K_{t}\right)-\frac{1}{p}\log\left(\frac{1}{|\mu|} \int_{S^{n-1}}h^p_{t}(u) d \mu(u)\right), \quad t \in(-\delta, \delta).
$$

Since $K_0$ is a maximizer of problem (\ref{4.3}) and $h_{K_{t}} \leq h_{t}$, we have
$\phi(0)=J_{m,p,\mu}\left(K_0\right) \geq J_{m,p,\mu}\left(K_{t}\right) \geq \phi(t)$
for all $t \in(-\delta, \delta)$. By Theorem \ref{2.1}, $\phi(t)$ is differentiable at $t=0$.
Therefore, from
$\phi'(0)=0$ and Theorem \ref{2.1}, we conclude that
\begin{equation}\label{pmPsi}
\left(m \widetilde{\Psi}_m\left(K_0\right)\right)^{-1} \int_{S^{n-1}} f(u) d \widetilde{C}_{p,m}^{\mathrm{a}}\left(K_0, u\right)=\left(\int_{S^{n-1}}h_{K_0}^p(u)d\mu(u)\right)^{-1}\int_{S^{n-1}} f(v) d \mu(u).
\end{equation}

Let $\bar{K}=c K_0$, where
\begin{equation}\label{cmnp}
c^{mn-p}=\left(m\widetilde{\Psi}_m\left(K_0\right)\right)^{-1}\int_{S^{n-1}}h_{K_0}^p(u)d\mu(u).
\end{equation}
From (\ref{pmPsi}), (\ref{cmnp}) and the $(m n-p)$-homogeneity of  $\widetilde{C}_{p,m}^{\mathrm{a}}$, we obtain
$$\int_{S^{n-1}} f(u) d \widetilde{C}_{p,m}^a(\bar{K}, u)=\int_{S^{n-1}} f(u) d \mu(u).$$
By the arbitrariness of $f \in C_e\left(S^{n-1}\right)$, we get the desired conclusion $\mu=\widetilde{C}_{p,m}^{\text {a }}(\bar{K}, \cdot)$.
\end{proof}



\begin{lemma}{\color{red}\cite[Lemma 4.2]{BLYZZ19}}\label{LL4.1}  Let $\varepsilon>0$. Suppose that $\left(e_{1, l}, \ldots, e_{n, l}\right)$, where $l=1,2, \ldots$, is a sequence of ordered orthonormal basis of $\mathbb{R}^n$ converging to the ordered orthonormal basis $\left(e_1, \ldots, e_n\right)$. Suppose also that $\left(a_{1, l}, \ldots, a_{n, l}\right)$ is a sequence of $n$-tuples satisfying
$$0<a_{1, l} \leq a_{2, l} \leq \cdots \leq a_{n, l}$$
for all $l$, and there exists an $\varepsilon>0$, such that $
a_{n, l}>\varepsilon$
for all $l$. For each $l=1,2, \ldots$, let
$$
Q_l=\left\{x \in \mathbb{R}^n:|x \cdot e_{1, l}|^2 / a_{1, l}^2+\cdots+|x \cdot e_{n, l}|^2 / a_{n, l}^2 \leq 1\right\}
$$
denote the ellipsoid generated by the $\left(e_{1, l}, \ldots, e_{n, l}\right)$ and $\left(a_{1, l}, \ldots, a_{n, l}\right)$. If $\mu$ is a Borel measure on $S^{n-1}$ that satisfies the strict subspace concentration inequality (\ref{E4.2}), then there exists $t_0, \delta_0, l_0>0$ and $c_{n,\varepsilon, t_0, \delta_0}$, independent of $l$, such that, for each $l>l_0$,
\begin{equation}\label{IEE}
E_{0,\mu}\left(Q_l\right) \leq-\frac{1}{n}\log \left(a_{1, l} \cdots a_{n-1, l}\right)+t_0 \log a_{1, l}+c_{n, \varepsilon, t_0, \delta_0}.
\end{equation}

\end{lemma}

The inequality (\ref{IEE}) in Lemma \ref{LL4.1} is still established if $E_{0,\mu}$ is replaced by $E_{p,\mu}$ when $p>0$. In fact, by Jessen's inequality, we have
\begin{eqnarray}\label{EEpmu}
E_{p,\mu}(Q_l)&=&-\frac{1}{p} \log\left(\frac{1}{|\mu|} \int_{S^{n-1}} h^p_{Q_l}(v) d \mu(v)\right)\\
&\leq&-\frac{1}{p} \left(\frac{1}{|\mu|} \int_{S^{n-1}}\log h^p_{Q_l}(v) d \mu(v)\right)\nonumber\\
&=&-\frac{1}{|\mu|} \int_{S^{n-1}}\log h_{Q_l}(v) d \mu(v)\nonumber\\
&=&E_{0,\mu}\left(Q_l\right).\nonumber
\end{eqnarray}
The proof of the following Lemma \ref{LL4.2} follows along the same lines as that of \cite[Theorem 5.2]{CLWX}, and thus we omit the proof.
\begin{lemma}\label{LL4.2} Suppose $p\geq0$, $m=1, \ldots, n-1$, and $\mu$ is an even finite Borel measure on $S^{n-1}$ with $|\mu|>0$. If $\mu$ satisfies the strict subspace concentration inequality (\ref{E4.2}),
then there exists $K_0 \in \mathcal{K}_{(e)}^n$ such that
\begin{equation}
J_{m,p,\mu}\left(K_0\right)=\sup \left\{J_{m,p,\mu}(K): K \in \mathcal{K}_{(e)}^n\right\}.\nonumber
\end{equation}
\end{lemma}

Theorem \ref{T1.2} on the solution of the Minkowski problem for measures $\widetilde{C}_{p,m}^a$ is a direct corollary of Lemma \ref{T4.1} and Lemma \ref{LL4.2}.

\bigskip

{\bf Funding}: Youjiang Lin is supported by NSFC 11971080, NSFC 12371137, Hebei Normal University Doctoral Research Start-up Fund L2025B50 and Jiangxi Provincial Natural Science Foundation 20232BAB201005. Yuchi Wu was supported, in part, by NSFC 12401067 and Science and Technology Commission of Shanghai Municipality 22DZ2229014.




\end{document}